\documentclass[a4paper,12pt]{amsart}
\usepackage{amsfonts}
\usepackage{amsmath, dsfont, amsthm, amsfonts, amssymb, color}
\usepackage{mathrsfs}
\usepackage{color}
\usepackage{stmaryrd}
\usepackage[utf8]{inputenc}
\usepackage[T1]{fontenc}
\usepackage[UKenglish]{babel}
\usepackage{verbatim}

\allowdisplaybreaks[4]
\usepackage{times}
\usepackage{latexsym}
\usepackage{booktabs}
\usepackage{mathtools}

\usepackage[colorlinks,pagebackref]{hyperref}
\usepackage{xcolor}

\newcommand{\xiang}{\color{blue}}

\theoremstyle{definition}

\newcommand{\scr}[1]{\mathscr #1}
\definecolor{wco}{rgb}{0.5,0.2,0.3}

\numberwithin{equation}{section} \theoremstyle{remark}

\allowdisplaybreaks
\def\R{\mathbb R}  \def\ff{\frac} \def\ss{\sqrt} \def\B{\mathbf
B}
 \def\kk{\kappa} 

\def\dd{\delta} \def\DD{\Delta} \def\vv{\varepsilon} 
\def\<{\langle} \def\>{\rangle}  
  \def\nn{\nabla} \def\pp{\partial} \def\E{\mathbb E}
\def\d{\text{\rm{d}}}  \def\aa{\alpha} 
  \def\si{\sigma} 
\def\beg{\begin} \def\beq{\begin{equation}}  \def\F{\scr F}
 
\def\e{\text{\rm{e}}}  \def\OO{\Omega}  
 \def\tt{\tilde} 
 \def\P{\mathbb P} 
\def\C{\scr C}           
   
\def\Z{\mathbb Z}  
\def\L{\scr L}

 \def\B{\scr B}  
\def\to{\rightarrow}

\def\3{\triangle}
\def\8{\infty}\def\1{\lesssim}
 
 \def\R{\mathbb R}  \def\ff{\frac} \def\ss{\sqrt} \def\B{\mathbf
B} \def\W{\mathbb W}
 \def\kk{\kappa} 

\def\dd{\delta} \def\DD{\Delta} \def\vv{\varepsilon} 
\def\<{\langle} \def\>{\rangle}  
  \def\nn{\nabla} \def\pp{\partial} \def\E{\mathbb E}
\def\d{\text{\rm{d}}}  \def\aa{\alpha} 
  \def\si{\sigma} 
\def\beg{\begin} \def\beq{\begin{equation}}  \def\F{\scr F}
 
\def\e{\text{\rm{e}}}  \def\OO{\Omega}  
 \def\tt{\tilde} 
 \def\P{\mathbb P} 
\def\C{\scr C}           
   
\def\Z{\mathbb Z}  
\def\L{\scr L}

 \def\B{\scr B}  
\def\to{\rightarrow}

\def\L{\scr L}
\def\R{\mathbb R}  \def\ff{\frac} \def\ss{\sqrt} \def\B{\mathbf
B}
 \def\kk{\kappa} 
\def\dd{\delta} \def\DD{\Delta} \def\vv{\varepsilon} 
\def\<{\langle} \def\>{\rangle}  
  \def\nn{\nabla} \def\pp{\partial} \def\E{\mathbb E}
\def\d{\text{\rm{d}}}  \def\aa{\alpha} 
  \def\si{\sigma} 
\def\beg{\begin} \def\beq{\begin{equation}}  \def\F{\scr F}
 
\def\e{\text{\rm{e}}}  \def\OO{\Omega}  
 \def\tt{\tilde} 
 \def\P{\mathbb P} 
\def\C{\scr C}           
   
\def\Z{\mathbb Z}  
\def\L{\scr L}

 \def\B{\scr B}  
\def\to{\rightarrow}
\def\8{\infty}\def\beq{\begin{equation}}

\begin{document}

\title{{\bf   Unadjusted Langevin Algorithms for   SDEs with H\"older Drift}
 }
\thanks{F. Y. Wang is supported in
 part by   the National Key R\&D Program of China (2022YFA1006000, 2020YFA0712900) and   NNSFC (11921001). L. Xu is supported by NNSFC grant (No. 12071499) and University of Macau grant MYRG2018- 00133-FST}

\author{Xiang Li}
\address[X. Li]{1. Department of Mathematics, Faculty of Science and Technology, University of Macau, Macau S.A.R., China. 2. Zhuhai UM Science \& Technology Research Institute, Zhuhai, China}
\email{yc07904@um.edu.mo}

\author{Feng-Yu Wang}
\address[F. Y. Wang]{Center for Applied Mathematics, Tianjin
University, Tianjin 300072, China}
\email{wangfy@tju.edu.cn}

\author{Lihu Xu}
\address[L. Xu]{1. Department of Mathematics, Faculty of Science and Technology, University of Macau, Macau S.A.R., China. 2. Zhuhai UM Science \& Technology Research Institute, Zhuhai, China}
\email{lihuxu@um.edu.mo}


\maketitle

\begin{abstract}
Consider the following stochastic differential equation for $(X_t)_{t\ge 0}$ on $\R^d$ and its  Euler-Maruyama (EM)  approximation $(Y_{t_n})_{n\in \mathbb Z^+}$:
\beg{align*} &\d X_t=b( X_t) \d t+\sigma(X_t)\d B_t, \\
& Y_{t_{n+1}}=Y_{t_{n}}+\eta_{n+1} b(Y_{t_{n}})+\sigma(Y_{t_{n}})\left(B_{t_{n+1}}-B_{t_{n}}\right), \end{align*}
where $b:\mathbb{R}^d  \mapsto \mathbb{R}^d,\ \ \sigma: \R^d \to\mathbb{R}^{d \times d}$ are measurable,  $B_t$ is the $d$-dimensional Brownian motion, $t_0:=0,t_{n}:=\sum_{k=1}^{n} \eta_{k}$ for  constants $\eta_k>0$ satisfying $\lim_{k \rightarrow \infty} \eta_k=0$ and $\sum_{k=1}^\infty\eta_k   =\infty$. Under (partial) dissipation conditions ensuring the ergodicity, we obtain explicit convergence rates of  $\W_p(\mathscr{L}(Y_{t_n}), \mathscr{L}(X_{t_n}))+\W_p(\mathscr{L}(Y_{t_n}), \mu)\to 0$ as $n\to\infty$, where $\W_p$ is the $L^p$-Wasserstein distance for certain $p\in [0,\infty)$,   $\L(\xi)$ is the distribution of random variable $\xi$, and $\mu$ is the unique invariant probability measure of $(X_t)_{t \ge 0}$.  Comparing with the existing results where $b$ is at least $C^2$-smooth, our estimates apply to H\"older continuous drift and can be sharp in several specific situations.

 \end{abstract}

\setcounter{tocdepth}{2}

\tableofcontents

\section{Introduction}

We consider the following time homogenous stochastic differential equation(SDE) on $\R^d$:
\begin{align} \label{SDE}
    \d X_t=b ( X_t) \d t+\sigma  (X_t)\d B_t, \ \ t\ge 0,
\end{align}
where $b:\mathbb{R}^d \mapsto \mathbb{R}^d,\ \ \sigma: \R^d \to\mathbb{R}^{d \times d}$ are measurable, and $B_t$ is the $d$-dimensional Brownian motion under a probability base $(\OO,\F,\{\F_t\}_{t\ge 0}, \P)$. Under  conditions allowing singular drift,  the well-posedness, regularity estimates and exponential ergodicity have been derived  for \eqref{SDE}, see for instance     \cite{G-M,R23,X-Z,Z05,W23,XZ} and references therein. 

 In recent years,   stochastic  algorithms arising from  statistics and machine learning
 have been intensively developed to simulate the invariant probability measure for a stochastic system,  where a   practical algorithm  is the   unadjusted Langevin algorithm (ULA)   using the Euler-Maruyama(EM) scheme of the associated SDE. See for instance \cite{E-J-H,D17,D-M17,M-M-N,D-M19,P-P} for the background of the study.

 For a sequence of step sizes $\{ \eta_{k} \}_{k\ge 1}$, the EM Scheme of $(\ref{SDE})$ is defined by the following induction
\begin{align} \label{EM-1}
    Y_{t_{k+1}}=Y_{t_{k}}+\eta_{k+1} b(Y_{t_{k}})+\sigma(Y_{t_{k}})\left(B_{t_{k+1}}-B_{t_{k}}\right), \qquad \ \ Y_{t_0}=Y_0=X_0,\  k \in \mathbb Z^+,\ 
\end{align}
with $Y_{0}=X_{0}=x$, where $t_{k}:=\sum_{i=1}^{k} \eta_{i}$ and $B_{t_{k+1}}-B_{t_{k}}$ can be identified as $\sqrt{\eta_{k+1}} \zeta_k$  for   i.i.d. $d$-dimensional standard normal distributed random variables 
$\{\zeta_k\}_{k\ge 1}.$   The associated continuous time Euler Scheme is defined by
\begin{align} \label{EM}
    Y_{t}=Y_{t_{k}}+(t-t_{k})b(Y_{t_{k}})+\sigma(Y_{t_{k}})\left( B_{t}-B_{t_{k}}\right), \qquad t\in [t_{k},t_{k+1}), k\ge 0, Y_0=X_0.
\end{align}  \par

In addition to the aforementioned ULA algorithm, the equation (\ref{EM-1}) can also be interpreted as a noisy gradient descent(GD) algorithm as $b(x)=-\nabla V(x)$ with $V(x)$ being a loss function in a certain optimization problem. Motivated by these two algorithms and their variants, the Euler-Maruyama(EM) scheme of SDEs have been extensively investigated under different assumptions and settings, see for instance \cite{A-J-K, D17, D-M17, D-M19, M-F-W-B, R-T}.  

In this paper, we investigate the convergence rate of 
$$\W_p(\L(X_{t_n}), \L(Y_{t_n}))+ \W_p(\L(Y_{t_n}),\mu)\to 0\ \text{as}\   n\to\infty,$$ where $p\in [0,\infty)$, $\mu$ is the unique invariant probability measure of $(X_t)_{t \ge 0}$, $\L(\xi)$ is the distribution of a random variable $\xi$, and 
for $\C(\mu_1,\mu_2)$ being the class of couplings of probability measures $\mu_1,\mu_2$ on $\R^d$, 
 $$\W_p(\mu_1,\mu_2):=\inf_{\pi\in \C(\mu_1,\mu_2)} \bigg(\int_{\R^d\times\R^d} 1_{\{x\ne y\}}|x-y|^p\pi(\d x,\d y)\bigg)^{\ff 1 {p\lor 1}}$$ 
 is the $L^p$-Wasserstein distance.  By Kantorovich's dual formula, for any $p\in [0,1]$ we have
$$\W_p(\mu_1,\mu_2)=\sup_{f\in\B_b, [f]_p\le 1} |\mu_1(f)-\mu_2(f)|,$$
where $\B_b$ is the set of all bounded measurable functions on $\R^d$, and
$$[f]_p:=\sup_{x\ne y} \ff{|f(x)-f(y)|}{|x-y|^p}$$
is the $p$-order modulus of continuity for $f$. In particular,
$$2\W_0(\mu_1,\mu_2)=  \|\mu_1-\mu_2\|_{TV}:=\sup_{|f|\le 1} |\mu_1(f)-\mu_2(f)|$$
  is   the total variation distance. \par

  It is worthy to point out that 
 for $b(x)=-\nabla V(x)$   and $\si=I_d$ (the $d\times d$ identity matrix),   where   $V$ is a strictly convex function with bounded $\nn^2 V,$ 
 the convergence rate under   the Kullback-Leibler divergence (i.e. the relative entropy) has been studied in  \cite{D17, A-J-K, M-F-W-B}.

By the ergodicity, we have $\L(X_{t_n}) \to \mu$ as $t_n\to\infty$, where $\L(X_{t_n})$ is the distribution of $X_{t_n}$ and  $\mu$ is the unique invariant probability measure. So, it is essential to assume $t_n\to\infty$ as $n\to\infty$, i.e.
$\sum_{k=1}^\infty\eta_k=\infty.$ On the other hand, to approximate $\L(X_{t_n})$ by $\L(Y_{t_n})$, the step size $\eta_n$ should decay to $0$   as $n\to\infty$. 
Thus, throughout the paper  we make the following assumption.

\  

\begin{enumerate} \item[$(A0)$]
    The step size sequence $\{\eta_{k}>0\}_{k\ge 1}$ is non-increasing and satisfies
    \begin{align*}
        \lim_{k\to \infty} \eta_{k}=0, \qquad \sum_{k= 1}^\infty\eta_{k}=\infty. 
    \end{align*}
    \end{enumerate}

\ \newline
  A  typical examples of step size $\{\eta_{k}\}_{k\ge 1}$ satisfying $(A0)$ is $\eta_{k}=\theta k^{-a}$ for some $\theta\in (0,\infty)$ and $a\in (0,1]$.

 Before moving on, let us  recall some existing estimates  on the convergence rate presented in the literature for $\eta_{k}=\frac{1}{k}$.

When $b=-\nn V$  for $V$ being strictly convex and $\nn^2 V$ being bounded,  Durmus and Moulines   \cite{D-M17,D-M19} showed that $\W_2(\L(Y_{t_{n}}), \mu)$ and $\W_0(\L(Y_{t_{n}}), \mu)$ are both bounded by
$O( d^{\frac{1}{2}}n^{-\frac{1}{2}})$,
and that their bounds can be improved to $O(dn^{-1})$ and  $O(d^{\frac{1}{2}}n^{-\frac{3}{4}})$ respectively if $V\in C^{3}$ has Lipschitz continuous first and second order derivatives.

 \par
Recently, Pagès and Panloup \cite{P-P} studied the non-asymptotic bounds related to the EM scheme (\ref{EM-1}), for which they replaced the strongly convex assumption on $V$ in \cite{D-M19} with a weaker assumption on the drift $b$ and allows  a multiplicative  noise  (i.e. $\sigma$ is non-constant).   They obtained an $O(n^{-1+\epsilon})$ upper bound with $\epsilon\in (0,1)$, 
for   $\W_0(\mathscr{L}(Y_{t_{n}}),\mu)$ under the conditions $b,\sigma\in C^{6}$, and an $O(n^{-1}\log(n))$ upper bound for $\W_1(\mathscr{L}(Y_{t_{n}}),\mu)$ under the conditions $b,\sigma \in C^{4}$. 
When $\sigma  $ is a constant   matrix, their upper bounds for both distances can be improved to $O(n^{-1})$ under the weaker condition that $b\in C^{3}$. Note that their results fail to provide an explicit dependence on the dimension $d$. \par

However, in many important applications such as lasso and bridge regressions, the drift $b$ is much more singular than those required in the above references. This motivates us to study the convergence of the unadjusted Langevin algorithm for SDE \eqref{SDE} in the case  $b\in C^{\alpha}$ with $\alpha\in (0,2]$ for partially filling this gap. If $b \in C^{\alpha}$ with $\alpha \in (0,1]$, it means that $b$ is $\alpha$-H\"{o}lder continuous, while if $b \in C^{\alpha}$ with $\alpha \in (1,2]$, it means that $\nabla b$ is $(\alpha-1)$-H\"{o}lder continuous.

To compare our main results with the above introduced estimates,   we simply let $\eta_{k}=\frac{1}{k}$  and briefly summarize our results as follows. 
 \begin{itemize}
     \item[1.] When $b$ is partially dissipative with $b\in C^\aa$ for some $\alpha\in (0,2]$, and $\sigma (x)=\sigma$ is constant, Theorem \ref{T2} implies    that 
     $$\W_1(\L(Y_{t_{n}}), \mu) + \W_0(\L(Y_{t_{n}}), \mu)\le  O(d^{\frac{1}{2}+1_{\{1<\alpha\le 2\}}}n^{-\frac{\alpha}{2}}).$$
 This estimate is also extended to $\W_p$ for $p\in (0,1)$.When $\alpha=2$, and $p=0,1$, the order of $n$ in our result matches the optimal order shown in \cite[Theorem 2.3]{P-P}. Note that our result is obtained under   weaker condition on  $b$ rather than $b \in C^3$ as in \cite{P-P}.
 \item[2.] When $b$ is uniformly dissipative and in $C^\aa$ for some $\alpha\in (0,1]$, Theorem \ref{T1} implies  
   $$\W_p(\L(Y_{t_{n}}), \L(X_{t_{n}}))  \le  O(d^{\frac{1}{2}}n^{-\frac{\alpha }{2}})$$ for any $p\in (1,\infty).$ 
   In particular, when $\alpha=1$ and $p=2$ this estimate has the same order
as that in \cite{D-M19}.
 \end{itemize}
\par
Let us briefly discuss our approach to proving the main results. In the uniform dissipation case,  we use a synchronuous coupling of the solutions of (\ref{SDE}) and (\ref{EM}) and compare their difference directly. However, this method is not applicable in the case of partial dissipation. Alternatively,  we apply a domino decomposition in \cite{P-P}, whose probability version was recently established in \cite{C-S-X} by a generalized Lindeberg principle and seems more powerful in applications. Comparing to \cite{P-P}, it seems to us that our one-step error estimation is a little more delicate, and thus allows us to achieve the optimal rate under a weaker assumption on $b$. More precisely, we establish an integration by parts formula for the unadjusted Langevin algorithm rather than directly estimate a Malliavin matrix, which helps us to get estimates explicitly depending on the dimension $d$.  Girsanov's theorem and Pinsker's inequality also play an important role in our proof.

\par
The paper is organized as follows.  When the coefficients only satisfy a  partial dissipation condition, we consider the additive noise where $\si (x)=\si$ does not depend on $x$ and $t$,
and   estimate $\W_p(\L(Y_{t_n}),\mu)$ in Section 2 for $p\in [0,1].$ In Section 3, we estimate   $\W_p(\L(Y_{t_n}),\L(X_{t_k}))$ for $p>1$ under the uniform dissipation condition, where the time in-homogenous  case  is also considered.

 Throughout   this paper, let $|\cdot|$ be the Euclidean norm,   let $\|\cdot\|_{{\rm op}}$ and $\|\cdot\|_{{\rm HS}}$ denote the operator norm and Hilbert-Schmidt norm for matrices respectively.
We will use $c=c(\cdots)$ to  denote a constant $c$ which only depends on the quantities   $``\cdots$".  The notations $C^k(\mathcal{X},\mathcal{Y})$ is used to denote the classes of continuous functions mapping from $\mathcal{X}$ to $\mathcal{Y}$ that have continuous first, second,\dots and $k$-th order partial derivatives and if the context of the function's space is clear, we will abbreviate the notation of the space as $C^k$.

\section{$\W_p$-estimate for $p\in [0,1]$: partial dissipation case   }

\subsection{Main result and an example}

\ \newline 
In this section, we study the $\W_p$-estimate for $p\in [0,1]$ under the following assumptions
\begin{enumerate} \item[$(A1)$]
    Let $\alpha\in (0,2]$.   $\si (x)=\si$ does not depend on $x $ and  is invertible, and there exist constants $K_1,K_2\in (0,\infty)$ such that
    $$\|\si\|_{\mathrm{op}}\lor \|\si^{-1}\|_{\mathrm{op}}\le K_1$$
    and $b$ is partially dissipative, i.e. 
    $$\< b(x)-b(y), x-y\>\le K_1-K_2|x-y|^2,\ \ \ x,y\in\mathbb R^d.$$ Moreover,
  \item[$\bullet$] when $\aa\le 1$,
 $$|b(x)-b(y)|\le K_1(|x-y|+|x-y|^\aa),\ \ \ x,y\in \R^d;$$
\item[$\bullet$]  when  $\aa\in (1,2]$, $b\in C^\aa(\R^d)$ and
    $$  \|\nn b\|_{\mathrm{op}}\le K_1,\ \       \|\nabla b(x)-\nabla b(y)\|_{\mathrm{op}}\le K_1   |x-y|^{\alpha-1},\ \
     \ \ x,y\in\mathbb R^d. $$
  \end{enumerate}
 This assumption implies the well-posedness of \eqref{SDE}, see \cite{R23,Zhang11} for the well-posedness of more singular SDEs.
 The following result     improves \cite[Theorem 2.2]{P-P}.  In particular,  when $\aa=2$ and $p>0,$ the order $n^{-1}$ in \eqref{W} is sharp for $\eta_k\sim \ff 1 k$ is sharp since it is reached by the Ornstein-Uhlenbeck process,  for which  $\W_0(\L(Y_{t_n}^x),\mu)$ can be computed explicitly (see \cite[Section 4.6]{P-P}).

\beg{thm} \label{T2} Assume $(A0)$ and $(A1)$ for some $\aa\in (0,2].$ Then there exist constants $c_1=c_1(K_1,K_2,\eta,\aa), c_2=c_2(K_1,K_2)\in (0,\infty)$ such that
the following statements hold for any $x\in\R^d$ and $n\ge 2$.
\beg{enumerate} \item[$(1)$] Let $\aa\in (0,1].$  We have
\beq\label{W1} \beg{split}&\W_{1}\left(\scr L(X_{t_{n}}^x), \L(Y_{t_{n}}^x)\right) +    \W_{1}\left(\mu, \L(Y_{t_{n}}^x)\right) \\
    & \le  c_1  d^{\ff 1 2}(1+|x|)   \sum_{k=1}^{n-1}    \e^{ -c_{2}(t_{n}-t_{k}) } \eta_{k}^{1+\ff \aa 2},\end{split}\end{equation}
\beq\label{W0}\beg{split}& \W_{0}\left(\scr L(X_{t_{n}}^x), \L(Y_{t_{n}}^x)\right) +    \W_{0}\left(\mu, \L(Y_{t_{n}}^x)\right) \\
    & \le  c_1  d^{\ff 1 2}(1+|x|) \left[\sum_{k=1}^{n-1}    \e^{ -c_{2}(t_{n}-t_{k}) } (t_n-t_k)^{-\ff 1 2} \eta_{k}^{1+\ff \aa 2}+\eta_n^{\ff {1  +(\aa\land 1)}2} \right].
    \end{split}\end{equation}
    In particular, if we choose $\eta_{k}=\frac{\theta}{k}$ for some constant $\theta>\ff{\alpha}{2c_2}$, then there exists a constant $c_1'=c_1'(K_1,K_2,\theta,\aa)\in (0,\infty)$ such that
\beq\label{W}  \W_{i}\left(\scr L(X_{t_{n}}^x), \L(Y_{t_{n}}^x)\right) +     \W_{i}\left(\mu, \L(Y_{t_{n}}^x)\right)
        \le  c_1'  d^{\ff 1 2} n^{-\ff \aa 2}, \ \ i=0,1.  \end{equation}
 \item[$(2)$] Let  $\aa\in (1,2].$        We have
 \beq\label{W1'} \beg{split}&\W_{1}\left(\scr L(X_{t_{n}}^x), \L(Y_{t_{n}}^x)\right) +    \W_{1}\left(\mu, \L(Y_{t_{n}}^x)\right) \\
     &\le  c_1  d^{\frac{3}{2}} (1+|x|^2)   \sum_{k=1}^{n-1}    \e^{ -c_{2}(t_{n}-t_{k}) } \eta_{k}^{1+\ff \aa 2},\end{split}\end{equation}
\beq\label{W0'} \beg{split}&\W_{0}\left(\scr L(X_{t_{n}}^x), \L(Y_{t_{n}}^x)\right) +    \W_{0}\left(\mu, \L(Y_{t_{n}}^x)\right) \\
    & \le  c_1  d^{\frac{3}{2}} (1+|x|^2) \left[\sum_{k=1}^{n-1}    \e^{ -c_{2}(t_{n}-t_{k}) } (t_n-t_k)^{-\frac{1}{2}} \eta_{k}^{1+\ff \aa 2} +\eta_n^{\ff {1  +(\aa\land 1)}2}\right].\end{split}\end{equation}
    In particular, if we choose $\eta_{k}=\frac{\theta}{k}$ for some constant $\theta>\ff{\alpha}{2c_2}$, then there exists a constant $c_1'=c_1'(K_1,K_2,\theta,\aa)\in (0,\infty)$ such that
    \begin{align} \label{W'}
        \W_{i}\left(\scr L(X_{t_{n}}^x), \L(Y_{t_{n}}^x)\right) +     \W_{i}\left(\mu, \L(Y_{t_{n}}^x)\right)
        \le  c_1'  d^{\ff 3 2} n^{-\ff \aa 2}, \ \ i=0,1.
    \end{align}
    \item[$(3)$] For any $p\in (0,1)$,
    \begin{align*}
        &\W_{p}\left(\scr L(X_{t_{n}}^x), \L(Y_{t_{n}}^x)\right) +     \W_{p}\left(\mu, \L(Y_{t_{n}}^x)\right)  \\
        &\le c_{1}d^{1+\frac{p}{2}+1_{\{1<\alpha\le 2\}}} \left(1+|x|^{1+1_{\{1<\alpha\le 2\}}} \right)\left[\sum_{k=1}^{n-1}    \e^{ -c_{2}(t_{n}-t_{k}) } (t_n-t_k)^{-\frac{1}{2}} \eta_{k}^{1+\ff \aa 2} +\eta_n^{\ff {1  +(\aa\land 1)}2}\right].
    \end{align*}
  \end{enumerate}
  \end{thm}

\textbf{Example 2.1.} Let us consider using the ULA algorithm for sampling from a $d$-dimensional distribution with density function $\frac{1}{Z}e^{-\frac{1}{2}|x|^{2}+\frac{1}{4}|x|^{\alpha+1}}, \alpha\in (0,1]$ where $Z$ is the normalization constant. In this case, the iterative formula for the ULA algorithm with decreasing step size $\eta_{k},k\ge 1$ is as follow
\begin{align*}
    Y_{k+1}=Y_{k}+\eta_{k+1} \left(-Y_{k}+\frac{\alpha+1}{4}\left|Y_{k} \right|^{\alpha-1}Y_{k} \right)+\sqrt{2\eta_{k+1}} \zeta_{k+1}, \quad k\ge 0,
\end{align*}
where $\zeta_{1},\dots, \zeta_{k},\dots$ are i.i.d. $d$-dimensional standard normal distributed random variables, and the corresponding SDE should be
\begin{align*}
    \mathrm{d}X_{t}=\left(-X_{t}+\frac{\alpha+1}{4}|X_{t}|^{\alpha-1}X_{t}\right)\mathrm{d}t+\sqrt{2}\mathrm{d}B_{t}, \quad t\ge 0.
\end{align*}
It can be easily verified that in this case, the drift coefficient $b(x)=-x+\frac{\alpha+1}{4}|x|^{\alpha-1}x$ is not uniformly dissipative (in other words, $V$ is non-convex for $V$ satisfying  $b=-\nabla V$), and non-differentiable near the origin if $\aa<1$. However, this example satisfies our assumptions in $(A1)$. Hence, according to Theorem \ref{T2}, if $\eta_{k}$ is chosen to be $\frac{\theta}{k}$ for some suitable $\theta>0$, then the approximation error of this algorithm under $\W_1+\W_0$ has convergence rate   $O(d^{\frac{1}{2}}n^{-\frac{\alpha}{2}})$.\\
  
  \subsection{Proof of Theorem \ref{T2}} 

\ \newline
We first present the following moment estimates on $Y_t$ and $X_t$ which are crucial for the proofs.

\beg{lem}\label{L2.2}  Assume that $(A0)$ holds, $\eqref{SDE}$ is well-posed and there exist positive constants $\kk_1,\kk_2$ such that
\beq\label{ZZ} \<x, b(x)\>\le \kk_1-\kk_2|x|^2,\ \ |b(x)|\le \kk_2 (1+|x|),\ \  \|\si\|_{\mathrm{op}}\le \kk_2,\ \ \forall x\in\R^d.\end{equation}
Then for any $p\in (0,\infty)$,  there exists a constant $\kk=\kk(\kk_1,\kk_2,p)\in (0,\infty)$ such that
\beq\label{Z1}\sup_{t\ge 0}\E |X_t^x|^p \le \kk (d^{\ff p 2}+|x|^p),\ \ \forall x\in\R^d.\end{equation}
If moreover
\beq\label{ZZ'} |b(x)-b(y)|\le K(|x-y|+|x-y|^\aa),\ \ \forall x,y\in \R^d \end{equation}
holds for some constant $K\in (0,\infty)$ and $\alpha\in (0,1]$, then there exists $\kk'=\kk'(\kk_1,\kk_2,K,\eta,\aa,p)\in (0,\infty)$ such that {\xiang 
\beq\label{Z2}\sup_{t\ge 0}\E |Y_t^x|^p \le \kk' (d^{\ff p 2}+|x|^p),\ \ \forall x\in\R^d.\end{equation}}
\end{lem}

\begin{proof} By Jensen's inequality, we only need to prove for $p\ge 2.$

 (1) Proof of \eqref{Z1}. All constants $c_i$ in this step  depend only on $\kk_1,\kk_2$ and $p$.
  By \eqref{ZZ}, $\|\si\|_{\mathrm{HS}}^2\le d \|\si\|_{\mathrm{op}}^2$  and Young's inequality, there exist constants $c_1,c_2\in (0,\infty)$ such that, for any $x,y\in \mathbb{R}^{d}$
 \beq\label{ZZ1}\beg{split}& p |x|^{p-2} \<x,b(x)\>+ \ff 1 2 p(p-1)|x|^{p-2}\|\si\|_{\mathrm{HS}}^2\\
 &\le -\kappa_{2}p|x|^{p}+\kappa_{1}p|x|^{p-2}+\frac{1}{2}\kappa_{2}p(p-1)d |x|^{p-2}\\
 &\le c_1 d^{\ff p 2 }-c_2|x|^p.\end{split}\end{equation}
 So, by It\^o's formula and $\|\sigma\|_{\mathrm{op}}\le \|\sigma\|_{\mathrm{HS}}$,
 \begin{align*}
     &\d |X_t^x|^p-\mathrm{d}M_{t}\\
     &=\bigg{[}p|X_{t}^{x}|^{p-2}\langle X_{t}^{x},b(X_{t}^{x})\rangle+\frac{1}{2}p|X_{t}^{x}|^{p-2}\|\sigma\|^{2}_{\mathrm{HS}}+\frac{1}{2}p(p-2)\left|X_{t}^{x} \right|^{p-4}\left| \sigma X_{t}^{x}\right|^{2}\bigg{]}\mathrm{d}t\\
     &\le \left[p|X_{t}^{x}|^{p-2}\langle X_{t}^{x},b(X_{t}^{x})\rangle+\frac{1}{2}p(p-1)|X_{t}^{x}|^{p-2}\|\sigma\|^{2}_{\mathrm{HS}}\right]\mathrm{d}t\\
     &\le \big(c_1d^{\ff p 2}-c_2|X_t^x|^p)\d t
 \end{align*}
 holds for some martingale $M_t$. Then it follows from Gronwall's inequality that
$$ \E|X_t^x|^p\le |x|^p\e^{-c_2 t}+c_1d^{\ff p 2} \int_0^t \e^{-c_2(t-s)} \d s\le \ff {c_1d^{\ff p 2}}{c_2}  + |x|^p,$$
so that \eqref{Z1} follows.

(2) Proof of \eqref{Z2}. All constants $c_i$ in this step  depend only on $\kk_1,\kk_2,K,\eta,\aa$ and $p$.
For any $k\ge 1$, by \eqref{ZZ}, \eqref{ZZ'} and \eqref{ZZ1},   and noting that 
$$  | Y_{t}^{x}-Y_{t_{k-1}}^{x} |+  | Y_{t}^{x}-Y_{t_{k-1}}^{x} |^{\aa}\le c(1+  | Y_{t}^{x}-Y_{t_{k-1}}^{x} |)$$
holds for some constant $c>0,$  we find a martingale $M_t$ such that for $t\in [t_{k-1},t_{k}]$
    \begin{align*}
       &\d |Y_t^x|^p -\d M_t
       \le \Big(p|Y_t^x|^{p-2}\<Y_t^x, b(Y^x_{t_{k-1}})\> +\ff 1 2 p(p-1)|Y_t^x|^{p-2}\|\si\|_{\mathrm{HS}}^2\Big)\d t\\
       &\le \Big(p|Y_t^x|^{p-2}\<Y_t^x, b(Y^x_{t})\> +\ff 1 2 p(p-1)|Y_t^x|^{p-2}\|\si\|_{\mathrm{HS}}^2\Big)\d t\\
       &\quad +p|Y_t^x|^{p-1} |b(Y_t^x)-b(Y^x_{t_{k-1}})| \d t  \\
       &\le \left[-\kappa_{2}p\left|Y_{t}^{x} \right|^{p}+\frac{1}{2}\kappa_{2}p(p-1)d\left| Y_{t}^{x}\right|^{p-2}+Kp\left| Y_{t}^{x}\right|^{p-1}\left( \left| Y_{t}^{x}-Y_{t_{k-1}}^{x}\right|+1\right)\right]\mathrm{d}t\\
       &\le \Big(-\ff{c_2}2 |Y_t^x|^p+   c_3 \big\{d^{\frac{p}{2}} + |Y_t^x-Y_{t_{k-1}}^x|^p \big\}\Big)\d t,
       \end{align*}
where the last  step follows from the Young's inequality. This implies
 \beq\label{Z2*}\beg{split}&\E|Y_{t_k}^x|^p\le \e^{-\ff 1 2 c_2 \eta_k} \E|Y_{t_{k-1}}^x|^p
   + c_3 d^{\frac{p}{2}} \eta_k \\
   &\quad + 2c_3 \int_{t_{k-1}}^{t_k}
 \E[|Y_t^x-Y_{t_{k-1}}^x|^{p}]\d t,\ \ k\ge 1.\end{split}\end{equation}
 By the boundedness of $\si$ and the linear growth of $b$, we find a constant $c_4>0$ such that
  \beg{align*} &\E|Y_t^x-Y_{t_{k-1}}^x|^p =\E\big|(t-t_{k-1}) b(Y_{t_{k-1}}^x) +\si (B_t-B_{t_{k-1}})\big|^p    \\
 &\le c_4 \left[\eta_k^{p} (1+ \E|Y_{t_{k-1}}^x|^p)+ \eta_{k}^{\frac{p}{2}}d^{\frac{p}{2}}\right],\ \ \ t\in [t_{k-1}, t_k], k\ge 1.\end{align*}
Since $\{\eta_{k}\}_{k\ge 1}$ is non-increasiong, combining this with \eqref{Z2*}, we can find a constant $c_5>0$
 such that
\beq\label{Z3} \begin{split}
    \E|Y_{t_k}^x|^p &\le \Big(\e^{-\ff 1 2 c_2 \eta_k}+ 2c_{3}c_{4} \eta_k^{1+p}\  \Big) \E|Y_{t_{k-1}}^x|^p
+ c_{3} d^{\frac{p}{2}}\eta_k+2c_{3}c_{4}\left[\eta_{k}^{1+p}+\eta_{k}^{1+\frac{p}{2}}d^{\frac{p}{2}} \right] \\
    &\le \Big(\e^{-\ff 1 2 c_2 \eta_k}+ c_5 \eta_k^{1+p}\  \Big) \E|Y_{t_{k-1}}^x|^p
+ c_5 d^{\frac{p}{2}}\eta_k,\ \ k\ge 1.
\end{split}
\end{equation}
   Since $\eta_k\to 0$ as $k\to\infty$, we can find a $k_0\ge 1$ such that
 \beq\label{Z4}  0\le \e^{-\ff 1 2 c_2 \eta_k}+ c_5 \eta_k^{1+p} \le 1- \ff 1 4 c_2\eta_k,\ \ k\ge k_0.\end{equation} Consequently,
    \begin{align*}
      c_6:=\sup_{k\ge 1}  \prod_{i=1}^{k} \Big(\e^{-\ff 1 2 c_2 \eta_j}+ c_5 \eta_j^{1+p}\Big) <\infty.
    \end{align*}
   So, iterating the estimate \eqref{Z3}  yields
  \beq\label{Z5}\E|Y_{t_{k}}|^{p}\le
        c_6|x|^{p}+c_5 d^{\frac{p}{2}} \left[\eta_{k}+\sum_{i=1}^{k-1} \eta_{i}\prod_{j=i+1}^{k}\Big(\e^{-\ff 1 2 c_2 \eta_j}+ c_5 \eta_{j}^{1+p}\Big)\right].\end{equation}
   By \eqref{Z4} we obtain
    \begin{align*}
    &\sum_{i=1}^{k-1}\eta_{i}\prod_{j=i+1}^{k}\Big(\e^{-\ff 1 2 c_2 \eta_j}+ c_5 \eta_j^{1+p}\Big)
    \le c_6 \sum_{i=1}^{k_0-1}\eta_{i}+\sum_{i=k_0}^{k-1}\eta_{i}\prod_{j=i+1}^{k}\left(1-\ff 1 4 c_2\eta_j \right) \\
    &= c_6\sum_{i=1}^{k_0-1}\eta_{i}+\frac{4}{c_2}\sum_{i=k_0}^{k-1}\left[1-(1-\ff 1 4 c_2 \eta_{i})\right]\prod_{j=i+1}^{k}\left(1-\ff 1  4 c_2\eta_{j}\right) \\
    &=c_6\sum_{i=1}^{k_0-1}\eta_{i}+\ff 4 {c_2 } \sum_{i=k_0}^{k-1}\left[\prod_{j=i+1}^{k}\left(1-\ff 1 4 c_2\eta_{j}\right)-\prod_{j=i}^{k}\left(1-\ff 1 4 c_2 \eta_{j}\right)\right] \\
    &\le c_6\sum_{i=1}^{k_0-1}\eta_{i}+\frac{4}{c_2}<\infty,\ \ \forall k\ge k_0.
    \end{align*}
   Combining this with   \eqref{Z5},  we prove \eqref{Z2} for some constant $\kk'$ and $k\ge k_0$, while for $k\le k_0$ the estimate follows from   \eqref{Z3}.
\end{proof}

\beg{lem}\label{C2.3}
Assume that the conditions in Lemma $\ref{L2.2}$ hold. Then, for any $p>0$ there exists a constant  $\kappa=\kappa(\kk_1,\kk_2,p)\in (0,\infty)$ such that
$$  \E |X_t^x-x|^p \le \kk (d^{\ff p 2}+|x|^{p})(1\land t)^{\ff p 2 \land 1},\ \ x\in\R^d,\ \ t\ge 0.$$
If moreover $\eqref{ZZ'}$ holds,  then there exists $\kk'=\kk'(\kk_1,\kk_2,\eta,K,p,\aa)\in (0,\infty)$ such that {
\begin{align*}
    \E |Y_t^x-Y_{t_{k-1}}^x|^p \le \kk' (d^{\ff p 2}+|x|^p)\eta_k^{\ff p 2},\ \ x\in\R^d,\ k\ge 1, \ t\in [t_{k-1},t_k].
\end{align*} }
\end{lem}

\begin{proof}
    Again, by Jensen's inequality, we only need to prove for $p\ge 2.$\par
    (1) Applying It\^o's formula, (\ref{ZZ1}) and Young's inequality gives us
    \begin{align*}
        &\mathrm{d}\left|X_{t}^{x}-x \right|^{p}-\mathrm{d}M_{t}\\
        &\le \left[p|X_{t}^{x}-x|^{p-2}\langle X_{t}^{x}-x,b(X_{t}^{x})\rangle+\frac{1}{2}p(p-1)|X_{t}^{x}-x|^{p-2}\|\sigma\|^{2}_{\mathrm{HS}}\right]\mathrm{d}t\\
        &\le \left(c_{1}d^{\frac{p}{2}}-c_{2}\left| X_{t}^{x}-x\right|^{p}+p |X_{t}^{x}-x|^{p-2}\langle X_{t}^{x}-x,b(X_{t}^{x})-b(X_{t}^{x}-x)\rangle\right)\mathrm{d}t\\
        &\le \left(c_{1}d^{\frac{p}{2}}-c_{2}\left| X_{t}^{x}-x\right|^{p}+2K_{1}p |X_{t}^{x}-x|^{p-1}(1+|x|)\right)\mathrm{d}t\\
        &\le \left(c_{1}'\left(d^{\frac{p}{2}}+|x|^{p} \right)-c_{2}'\left|X_{t}^{x}-x \right|^{p} \right)\mathrm{d}t
    \end{align*}
    where the third inequality is a consequence of (\ref{ZZ'}) and $M_{t}$ is the Martingale term. Since $X_{0}^{x}-x=0$, it follows from Gronwall's inequality that
    \begin{align*}
        \mathbb{E}\left| X_{t}^{x}-x\right|^{p}\le \frac{c_{1}'(d^{\frac{p}{2}}+|x|^{p})}{c_{2}'}\left( 1-e^{-c_{2}'t} \right)
        \le \kappa (d^{\frac{p}{2}}+|x|^{p}) \left(1 \land t \right), \qquad \forall x\in \mathbb{R}^{d},t\ge 0,
    \end{align*}
    for some $\kappa=\kappa(\kappa_{1},\kappa_{2},p)$. \par
   (2) For the second one, notice that for $t\in [t_{k},t_{k+1})$, the conditional distribution of $Y_{t}^{x}-Y_{t_{k}}^{x}$ given ${Y_{t_{k}}^{x}} $ is 
   $$Y_{t}^{x}-Y_{t_{k}}^{x}|_{{Y_{t_{k}}^{x}}}\sim \mathcal{N}\left(
 (t-t_{k})b(Y_{t_{k}}^{x}), (t-t_{k})\sigma\sigma^{T}\right).$$
 Hence, it follows from (\ref{ZZ}) and (\ref{Z2}) that for some constants $c_3'=c_3'(p)>0$ and $c_4'=c_4'(p,\kk_2)>0,$ 
 \begin{align*}
     \mathbb{E}|Y_{t}^{x}-Y_{t_{k}}^{x}|^{p}
     &\le c_3'\left[(t-t_{k})^{p}\mathbb{E}|b(Y_{t_{k}}^{x})|^{p}+(t-t_{k})^{\frac{p}{2}}\|\sigma\|_{\mathrm{HS}}^{p}\right]\\
     &\le c_3'\kappa_{2}^{p}\left[(t-t_{k})^{p}\left(1+\mathbb{E}|Y_{t_{k}}^{x}|\right)^{p}+(t-t_{k})^{\frac{p}{2}}d^{\frac{p}{2}}\right]\\
     &\le c_4' \left(d^{\frac{p}{2}}+|x|^{p} \right) (t-t_{k})^{\frac{p}{2}}.
 \end{align*}
 So the proof is complete.
\end{proof}

\ \newline 
To prove Theorem \ref{T2}, we consider the SDEs \eqref{SDE} and \eqref{EM} on each time interval $[t_k,t_{k+1})$ for $k\in \mathbb Z^+,$ where $t_0:=0.$ 
For any $x\in \R^d$ and $k\in\mathbb Z^+,$ let $ (Y_{t_k,t}(x))_{t\in [t_k,t_{k+1}]}$ solve the SDE  
\beg{align*} \d Y_{t_k,t}(x)= b(x)\d t+\si \d B_t,\ \ \ X_{t_k,t_k}(x)=Y_{t_k,t_k}(x)=x,  t\in [t_k, t_{k+1}]. \end{align*} 
Define 
  \beg{align*}  Q_{t_k, t_{k+1}}f(x):= \E[f(Y_{t_k,t_{k+1}}(x))],\ \  Q_{t_k,t_n} := Q_{t_k, t_{k+1}} Q_{t_{k+1}, t_{k+2}}\cdots Q_{t_{n-1}, t_n},\ n\ge k+1.\end{align*}   
Correspondingly, for any $s\ge 0$ and $x\in \R^d$, we  let $\{X_{s,t}(x)\}_{t\ge s} $ solve the SDE
$$\d X_{s,t}(x)= b(X_{s,t}(x))\d t+\si \d B_t,\ \ t\ge s, X_{s,s}(x)=x.$$ Then the Markov semigroup $P_t$ associated with \eqref{SDE} satisfies 
\beq\label{SP}   P_{t-s}f(x)= P_{s,t}f(x):=\E[f(X_{s,t}(x))],\ \  \ t\ge s\ge 0.\end{equation}   Let $Q_{0,0}=P_0$ be the identity operator. We have the domino decomposition  
$$P_{t_n}- Q_{0, t_n}= \sum_{k=1}^n Q_{0, t_{k-1}} (P_{t_{k-1},t_k} -Q_{t_{k-1},t_k})P_{t_k,t_n},\ \ n\in\mathbb N.$$
Combining this with 
$$\W_{p}\left(\scr L(X_{t_{n}}^x), \L(Y_{t_{n}}^x)\right)=\sup_{[f]_p\le 1}  |\E[f(X_{t_n}^x)- f (Y_{t_n}^x)] |   = \sup_{[f]_p\le 1}|P_{0,t_n} f(x)- Q_{0,t_n} f (x)|,\ \ n\ge 1, $$
we derive 
\beq\label{SM}  \W_{p}\left(\scr L(X_{t_{n}}^x), \L(Y_{t_{n}}^x)\right)= \sup_{[f]_p\le 1} \Big|\sum_{k=1}^n Q_{0, t_{k-1}} (P_{t_{k-1},t_k} -Q_{t_{k-1},t_k})P_{t_k,t_n}f(x)\Big|.
\end{equation}

To prove Theorem \ref{T2}  using this formula, we need the following derivative estimates on $P_{t}. $

 \beg{lem}\label{L3.1} Assume $(A0)$ and $(A1)$. There exist constants $\kk_1=\kk_1(K_1,K_2),\kk_2=\kk_2(K_1,K_2)\in (0,\infty)$ such that
\beq\label{GR1}  \|\nn  P_tf\|_\infty\le \kk_1 \e^{-\kk_2 t} \|\nn f\|_\infty,\ \ t>0,\end{equation}
\beq\label{GR2}  \|\nn  P_tf\|_\infty\le \kk_1 \e^{-\kk_2 t} t^{-\ff 1 2} \|f\|_\infty,\ \ t>0.\end{equation}
Moreover, when $\aa\ge 1$,
\beq\label{GR3}  \|\|\nn^2  P_tf\|_{\mathrm{op}}\|_\infty\le \kk_1    \e^{-\kk_2 t} \big(t^{-\ff 1 2 }+ \ss d\, 1_{\{\aa<2\}}\big) \|\nn f\|_\infty,\ \ t>0.\end{equation}
\end{lem}

 \beg{proof}  (a)   By \cite[Corollary 2.3]{E16}, $(A1)$ implies
\begin{align*}
\sup_{\|\nn f\|_\infty\le 1} |P_tf(x)-P_tf(y)|&= \W_1(\L(X_t^x), \L(X_t^y)) \\
& \le \kk_1  \e^{-\kk_2 t}|x-y|\ \  \ \  \ \  \ \  \ \  \ \ \ \  \forall \ \  t\ge 0, x,y\in\R^d,
\end{align*}
for some constants $\kk_1,\kk_2\in (0,\infty)$  depending only on $K_1,K_2$. Consequently, \eqref{GR1} holds.

Next, by \cite[Theorem 3.4]{P-W}, there exists constant $k_1 \in (0,\infty)$  depending only on $K_1,K_2$ such that
$$\|\nn_{v} P_tf\|_\infty\le  k_1   t^{-\ff 1 2} \|f\|_\infty |v|,\ \  t\in (0,1],\forall v\in \mathbb{R}^{d}.$$
Combining this with \eqref{GR1} and the semigroup property of $P_t$,  we prove \eqref{GR2} for some constants $\kk_1,\kk_2\in (0,\infty)$  depending only on $K_1,K_2$.

 (b) 
 When $\aa=2$, we have $\|\nn b\|_{\mathrm{op}}\lor \|\nn^2 b\|_{\mathrm{op}}\le K_1$. In general, we let $\|\nn^2 b\|_{\mathrm{op}}\le K_1'$ for  some constant $K_1'$ possibly   different from $K_1$.
 Then for any  $v,w\in\R^d$,
 $$ \phi_t(v):= \nn_v X_t^x :=\lim_{\vv\downarrow 0} \ff{X_t^{x+\vv v}-X_t^x}\vv,\ \ \varphi_{t}(v,w):= \lim_{\vv\downarrow 0} \ff{\nn_v X_t^{x+\vv w}- \nn_v X_t^x}{\vv}$$
 exist and solve the equations
 \beg{align*} &\pp_t \phi_t(v)= \nn_{\phi_t(v)}b(X_t^x),\ \ \ \phi_0(v)=v,\\
 &\pp_t \varphi_t(v,w)=  \nn_{\varphi_t(v,w)}b(X_t^x) + \nn_{\phi_t (w)}\nn_{\phi_t(v)} b (X_t^x),\ \ \varphi_0(v,w)=0.\end{align*}
 Consequently,
\beq\label{N1} \sup_{|v|\le 1} |\phi_t(v)|\le \e^{K_1},\ \  \sup_{|v|,|w|\le 1} |\varphi_t(v,w)|\le K_1'\e^{3K_1},\ \ \ t\in [0,1].\end{equation}
 Combining this with the Bismut–Elworthy–Li formula
 $$\nn_v P_t f(x)= \ff 1 t \E\bigg[f(X_t^x) \int_0^t \<\si^{-1} \phi_s(v), \d B_s\>\bigg],$$
 we derive
 $$\nn_w\nn_v P_t f(x)=  I_1+I_2,$$
 where
\beq\label{N0}
\beg{split}
&I_1:= \ff 1 t\E\bigg[\<\nn f(X_t^x), \phi_t(w)\>  \int_0^t \<\si^{-1} \phi_s(v), \d B_s\>\bigg],\\
&I_2:=     \ff 1 t\E\bigg[  f(X_t^x)   \int_0^t \<\si^{-1} \varphi_s(v,w), \d B_s\>\bigg].
 \end{split}
 \end{equation}
By \eqref{N1} and $\|\si^{-1}\|_{\mathrm{op}}\le K_1$, for any $x\in\R^d$ and $t\in (0,1],$
\beq\label{N01} \beg{split} &|I_1|\le \|\nn f\|_\infty K_1 \e^{2K_1} t^{-\ff 1 2}, \\
&|I_2|= \ff 1 t \bigg|\E\bigg[  \big\{f(X_t^x)  -P_tf(x)\big\} \int_0^t \<\si^{-1} \varphi_s(v,w), \d B_s\>\bigg] \bigg|\\
&\ \ \ \ \ \ \le \ff 1 t \big(\E|f(X_t^x)  -P_tf(x)|^2\big)^{\ff 1 2} K_1K_1' \e^{3 K_1}.\end{split} \end{equation}
Noting that \eqref{GR1} implies
\beq\label{PC} \beg{split}  \E|f(X_t^x)  -P_tf(x)|^2&= P_t f^2(x)-(P_tf(x))^2 \\
&=\int_0^t \ff{\d}{\d s} P_s (P_{t-s} f)^2(x)\d s\\
&=\int_{0}^{t} P_{s}\left(L\left( P_{t-s}f\right)^{2}-2P_{t-s}f\cdot LP_{t-s}f \right) (x) \mathrm{d}s\\
&= \int_0^t P_s |\si^* \nn P_{t-s}f|^2 (x)\d s\le \|\nn f\|_\infty^2(K_1\kk_1)^2t,\ \ t>0,\end{split}\end{equation}
where {$L:={\rm tr} \{\si\si^*\nn^2\}+ b\cdot\nn$} is the generator associated with (\ref{SDE}). We derived \eqref{GR3} for $t\in (0,1]$ for some larger constant $\kk_1=\kk_1(K_1,K_2)$  since $K_1'=K_1$ under $(A1)$. And for $t\in (1,\infty)$, the desired result follows from
\begin{align*}
    \|\|\nabla^{2}P_{t}f\|_{\mathrm{op}}\|_{\infty}&=\|\|\nabla^{2}P_{1}P_{t-1}f\|_{\mathrm{op}}\|_{\infty} \\
    &\le \kappa_{1}e^{-\kappa_{2}}\|\nabla P_{t-1}f\|_{\infty}\le \kappa_{1}e^{-\kappa_{2}t}\|\nabla f\|_{\infty},
\end{align*}
where the last inequality is a consequence of (\ref{GR1}).\par
Now, let $\aa\in [1,2).$ Let $\tt b(x)=\E[b(x +B_1)].$ Then
\beq\label{N} \beg{split}&|\nn_v \tt b(x)|\le \E[|\nn_v b(x+ B_1)|]\le K_1,\ \ |v|\le 1,\\
&|b(x)-\tt b(x)|\le K_1  \E|B_1|\le K_1 \sqrt{d}.\end{split}\end{equation}
By the Bismut formula,
$$\nn_v \tt b(x)= \E\big[b(x+B_1) \<v,B_1\>\big].$$
 For $w\in \R^d$ with $|w|\le 1$, $\|\nn b\|_{\mathrm{op}}\le K_1$ yields
\beg{align*}&\big|\nn_w\nn_v \tt b(x)\big|= \bigg|\lim_{\vv\downarrow 0} \ff 1 \vv \E \Big[\big\{b(x+ B_1 +\vv w) - b(x+ B_1 ) \big\} \<v,B_1\>\Big]\bigg| \le K_1 .\end{align*}
Let  $\tt P_t$ generated by $\tt L:={\rm tr} \{\si\si^*\nn^2\}+\tt b\cdot\nn$. By the result for $\alpha=2$ case, we can find a constant $c_1=c_1(K_1)\in (0,\infty)$ such that the semigroup satisfies
\begin{align} \label{tP}
    \sup_{x\in\R^d, |v|\lor |w|\le 1,\|\nn f\|_\infty\le 1 } |\nn_w\nn_v \tt P_tf(x)|\le c_1  \, t^{-\ff 1 2},\ \ t\in (0,1].
\end{align}
Combining this with \eqref{GR1}, \eqref{N}, $\|\nabla b\|_{\mathrm{op}}\le K_1$ and  the formula
$$P_tf=\tt P_tf+\int_0^t \tt P_s \nn_{b-\tt b}P_{t-s}f\d s,$$
 we derive that for any $x\in \mathbb{R}^{d},|v|\lor |w|\le 1,t\in (0,1]$ and $\|\nn f\|_\infty\le 1$,
 \begin{align*}
     \left|\nabla_{w}\nabla_{v} P_{t}f(x) \right|
     &=\left|\nabla_{w}\nabla_{v}\tt P_tf(x)+\int_0^t \nabla_{w}\nabla_{v}\tt P_s \nn_{b-\tt b}P_{t-s}f(x)\d s \right|\\
     &\le c_1   t^{-\ff 1 2}+c_1  \int_0^t  s^{-\ff 1 2}   \|\nabla \nabla_{b-\tt b}P_{t-s}f\|_{\infty}\d s.
 \end{align*}
By denoting $\xi_{t}:=\sup_{x\in\R^d, |v|\lor |w|\le 1,\|\nn f\|_\infty\le 1 } |\nn_w\nn_v  P_tf(x)|$, we then have
$$\xi_t
\le c_1   t^{-\ff 1 2}+c_1 K_1  \sqrt{d}\int_0^t  s^{-\ff 1 2}    \xi_{t-s}\d s,\ \ t\in (0,1].$$
Let $a_{t}=\sup_{s\in[0,t]}\{s^{\frac{1}{2}}\xi_{s}\}$ and notice that $\int_{0}^{1} s^{-\frac{1}{2}} \left(1-s \right)^{-\frac{1}{2}} \mathrm{d}s= B(1,1)=1$ with $B(\cdot,\cdot)$ being the beta function. It follows that
\begin{align*}
    a_{t}&\le c_{1} +c_{1}\sqrt{d}K_{1}t^{\frac{1}{2}} \left(\int_{0}^{t} s^{-\frac{1}{2}} (t-s)^{-\frac{1}{2}} \mathrm{d}s\right) a_{t} \\
    &
    =c_{1} +c_{1}\sqrt{d}K_{1}t^{\frac{1}{2}} B(1,1) a_{t}.
\end{align*}
Solving this inequality yields, for any $t\in \left(0,t_{0}/2\right)$ with $t_{0}=\frac{1}{(c_{1}K_{1})^{2}d}$,
\begin{align*}
    a_{t}\le \frac{c_{1}}{1-c_{1}\sqrt{d}K_{1}t^{\frac{1}{2}}}\le 2c_{1},
\end{align*}
which implies
\begin{align*}
    \|\|\nabla^{2}P_{t}f\|_{\mathrm{op}}\|_{\infty}=\xi_{t} \|\nabla f\|_{\infty}\le 2c_{1}t^{-\frac{1}{2}}\|\nabla f\|_{\infty}.
\end{align*}
Similar as before, for $t\in \left(t_{0}/2,\infty\right)$, the proof is finished by
\begin{align*}
    \|\|\nabla^{2}P_{t}f\|_{\mathrm{op}}\|_{\infty}&=\|\|\nabla^{2}P_{\frac{t_{0}}{4}}P_{t-\frac{t_{0}}{4}}f\|_{\mathrm{op}}\|_{\infty}\\
    &\le 4c_{1}t_{0}^{-\frac{1}{2}}\|\nabla P_{t-\frac{t_{0}}{4}}f\|_{\infty}\le \kappa_{1} \sqrt{d} t^{-\frac{1}{2}} e^{-\kappa_{2}t} \|\nabla f\|_{\infty},
\end{align*}
for some $\kappa_{1},\kappa_{2}\in (0,\infty)$.
 \end{proof}

\beg{lem}\label{L3.2} Assume $(A0)$ and $(A1)$. Then there exists  a constant  $\kk=\kk(K_1,K_2,\eta,\aa) $ such that the following statements hold for any $x\in\R^d$ and $k\in\mathbb N.$
\beg{enumerate} \item[$(1)$] When $\aa\in (0,1]$,
\beg{align*}
&\sup_{\|\nn f\|_\infty\le 1} \left| Q_{0,t_{k-1}} (P_{t_{k-1},t_k} -Q_{t_{k-1},t_k} ) P_{t_k,t}  f (x) \right|  \\
& \le  \kk d^{\ff 1 2}(1+|x|)\e^{-\kk_2(t-t_k)}  \eta_k^{1+\ff\aa 2},\ \  \ t>t_k.
 \end{align*}
\item[$(2)$] When $\aa\in (1,2]$,
 \beq\label{ES2}
 \beg{split}&\sup_{\|\nn f\|_\infty\le 1} \left| Q_{0,t_{k-1}} (P_{t_{k-1},t_k} -Q_{t_{k-1},t_k} ) P_{t_k,t}  f (x) \right|  \\
& \le  \kk d^{\frac{3}{2}} (1+|x|^2)\e^{-\kk_2(t-t_k)} \bigg( \eta_k^{1+\ff\aa 2}+ (t-t_k)^{-\ff 1 2}\eta_k^3 \bigg),\ \  \ t>t_k.\end{split}
\end{equation}

\end{enumerate}
\end{lem}

\beg{proof} (a) Let $\aa\in (0,1].$  By Lemma \ref{C2.3} and  and Lemma \ref{L2.2}, we can find a constant $k_1=k_1(K_1,K_2)\in (0,\infty)$ such that
\beq\label{AB} \beg{split}  &\ \ \ \ \ \E[ |X_{t_{k-1},t }(Y_{t_{k-1}}^x) - Y_{t_{k-1}}^x |+  |X_{t_{k-1},t }(Y_{t_{k-1}}^x) - Y_{t_{k-1}}^x|^\aa ]  \\
&= \E\Big\{ \big(\E[ |X_{t_{k-1},t_k }(z) - z |+  |X_{t_{k-1},t_k }(z) - z|^\aa ] \big)_{z=Y_{t_{k-1}}^x}\Big\}\\
&\le \kappa \left[ \left(d^{\frac{1}{2}}+\mathbb{E}\left| Y_{t_{k-1}}^{x}\right| \right)\eta_{k}^{\frac{1}{2}} + \left(d^{\frac{\alpha}{2}}+\mathbb{E}\left| Y_{t_{k-1}}^{x}\right|^{\alpha} \right)\eta_{k}^{\frac{\alpha}{2}}\right] \\
& \le k_1 d^{\frac{1}{2}}(1+|x|) \eta_k^{\ff\aa 2},\ \ \ t\in [t_{k-1},t_k].\end{split}
\end{equation}
Meanwhile, by (\ref{SDE}), (\ref{EM-1}), \eqref{GR1}  and  $(A1)$, we obtain
\beg{align*} &\ \ \ \ \ \sup_{\|\nn f\|_\infty\le 1} \left| Q_{0,t_{k-1}} (P_{t_{k-1},t_k} -Q_{t_{k-1},t_k} ) P_{t_k,t}  f (x) \right|   \\
&\le \|\nabla P_{t_{k},t}f\|_{\infty} \sup_{\|\nn g\|_\infty\le 1} \left| Q_{0,t_{k-1}} (P_{t_{k-1},t_k} -Q_{t_{k-1},t_k} )   g (x) \right|  \\
&\le \kk_1 \e^{-\kk_2(t-t_k)} \sup_{\|\nn g\|_\infty\le 1} \left| \E[ g(X_{t_{k-1},t_k }(Y_{t_{k-1}}^x) ) - g(Y_{t_{k-1},t_k }(Y_{t_{k-1}}^x) )]\right| \\
&\le  \kk_1 \e^{-\kk_2(t-t_k)}  \E \left| X_{t_{k-1},t_k }(Y_{t_{k-1}}^x) ) - Y_{t_{k-1},t_k }(Y_{t_{k-1}}^x) \right|  \\
&=   \kk_1 \e^{-\kk_2(t-t_k)} \int_{t_{k-1}}^{t_k} \E[ | b( X_{t_{k-1},t }(Y_{t_{k-1}}^x) ) - b(Y_{t_{k-1}}^x) | ]\d t\\
&\le \kk_1 K_1 \e^{-\kk_2(t-t_k)} \int_{t_{k-1}}^{t_k} \E \Big\{\big(\E[ |X_{t_{k-1},s}(z) - z |+  |X_{t_{k-1},s}(z) - z|^\aa ] \big)_{z=Y_{t_{k-1}}^x} \Big\}\d s.\end{align*}
Combining this with \eqref{AB} we prove the first assertion.

(b) Let $\aa\in (1,2].$
For any function $g$ on $\R^d$ with $\|\nn^i g\|_\infty<\infty, i=1,2,$    we have
\begin{align*}
    g(z)-g(y)
    &=\int_{0}^{1}\nabla_{z-y}g(y+r(z-y)) \mathrm{d}r\\
    &=\int_{0}^{1}\nabla_{z-y}g(y+r(z-y))-\nabla_{z-y}g(y) \mathrm{d}r+\left(\nabla_{z-y}g(y)-\nabla_{z-y}g(u)\right)+\nabla_{z-y}g(u)\\
    &=\int_{0}^{1}\int_{0}^{1}r_{1}\nabla_{z-y}\nabla_{z-y}g(y+r_{1}r_{2}(z-y)) \mathrm{d}r_{2}\mathrm{d}r_{1}\\
    &\quad +\int_{0}^{1}\nabla_{y-u}\nabla_{z-y}g(u+r(y-u))\mathrm{d}r+\nabla_{z-y}g(u),\ \ u,y,z\in\R^d.
\end{align*}
Let $\|\nn f\|_\infty\le 1$. We shall apply this formula for
$$g= P_{t_k,t} f,\ \ z= X_{t_{k-1},t_k}(Y_{t_{k-1}}^x),\  \ y= Y_{t_{k-1}, t_k}(Y_{t_{k-1}}^x),\ \ u=Y_{t_{k-1}}^x.$$
Let
\begin{align} \label{D}
    \DD_k:= X_{t_{k-1},t_k}(Y_{t_{k-1}}^x) - Y_{t_{k-1},t_k}(Y_{t_{k-1}}^x),\ \ \  \tt \DD_k:= Y_{t_k}^x-Y_{t_{k-1}}^x.
\end{align}
 By \eqref{SP} and noting that $Y_{t_{k-1},t_k}(Y_{t_k-1}^x)=Y_{t_k}^x,$ we deduce  from the above formula that
\beq\label{TT}  \beg{split}  &\  \ \ \ \ \big|Q_{0,t_{k-1} } (P_{t_{k-1},t_k} -Q_{t_{k-1},t_k} ) P_{t_k,t}  f (x)\big|\\
    &=\left| \mathbb{E}\left[P_{t_{k},t}f\left( X_{t_{k-1},t_{k}} (Y_{t_{k-1}}^x)\right) \right]-\mathbb{E}\left[P_{t_{k},t_{n}}f\left( Y_{t_{k-1},t_{k}} (Y_{t_{k-1}}^x)\right) \right] \right|\\
    &\le \left|\mathbb{E} \int_{0}^{1}\int_{0}^{1} r_{1}\nabla_{\Delta_{k} }\nabla_{\Delta_{k} }P_{t_{k},t} f\left(Y_{t_{k-1},t_{k}}  (Y_{t_{k-1}}^x)+r_{1}r_{2}\Delta_{k}  \right)\mathrm{d}r_{2}\mathrm{d}r_{1} \right|  \\
    & \ \ \ \ +\left|\mathbb{E} \int_{0}^{1} \nabla_{\tt \Delta_{k}}  \nabla_{\Delta_{k} }P_{t_{k},t}f \left(Y_{t_{k-1}}^x+r\left(Y_{t_{k}}^x-Y_{t_{k-1}}^x\right) \right) \mathrm{d}r \right|  \\
    &\ \ \ \  +\left|\mathbb{E}[\nabla_{\Delta_{k}}P_{t_{k},t}f (Y_{t_{k-1}}^x) ]\right|  \\
    &:=J_{1}+J_{2}+J_{3},  \end{split}\end{equation}
    To estimate these terms, we  need to bound $ \DD_k$ and $\tt \DD_k.$ For $t\in [t_{k-1},t_{k}]$, define
    \begin{align}
          \DD_{k,t}&:=X_{t_{k-1},t}(Y_{t_{k-1}}^x)-  Y_{t_{k-1},t}(Y_{t_{k-1}}^x) = X_{t_{k-1},t}(Y_{t_{k-1}}^x)-  Y_{t}^x,\label{KT}\\
         \tt \DD_{k,t}&:=Y_{t_{k-1},t}(Y_{t_{k-1}}^x)-  Y_{t_{k-1}}^x = Y_{t}^x-Y_{t_{k-1}}^{x} .\label{tKT}
    \end{align}
    By Lemma \ref{C2.3},  we can find a constant $C_0=C_0(K_1,K_2,\eta,\aa)\in (0,\infty)$ such that
\beq\label{D0} \E|\tt \DD_{k,t}|^4 \le C_0 d^2 (1+|x|^4) (t-t_{k-1})^2,\quad \forall t\in [t_{k-1},t_{k}].\end{equation}
We claim that there exist  constants $C_i=C_i(K_1,K_2,\eta,\aa)\in (0,\infty), i=1,2,$ such that
\beq\label{D1} \sup_{t\in [t_{k-1},t_k]} \E|\DD_{k,t}|^4  \le C_1 d^{2}  (1+|x|^4) \eta_k^{6},\end{equation}
\beq\label{D2}\E \big|\E[\DD_k|Y_{t_{k-1}}^x]\big|  \le C_2 d^{\ff \aa2} (1+|x|^{\aa})  \eta_k^{1+\ff\aa 2}.\end{equation}
Indeed,  by    \eqref{D0} and  that $(A1)$ for $\aa\in (1,2]$,
\beg{align*}
 \E|\DD_{k,t} |^4&=\E| X_{t_{k-1},t}(Y_{t_{k-1}}^x) - Y_{t_{k-1},t}(Y_{t_{k-1}}^x)|^4 \\
&\le \eta_k^3 \int_{t_{k-1}}^t \E|b(X_{t_{k-1},s}(Y_{t_{k-1}}^x))- b( Y_{t_{k-1}}^x) |^4 \d s \\
&\le 8 K_1^2   \eta_k^3 \int_{t_{k-1}}^t \E\big\{  | X_{t_{k-1},s}(Y_{t_{k-1}}^x)-  Y_{t_{k-1},s}(Y_{t_{k-1}}^x)|^4+    |Y_{t_{k-1},s}(Y_{t_{k-1}}^x)-   Y_{t_{k-1}}^x|^4 \big\}\d s \\
& \le 8 K_1^2   \eta_k^3  \int_{t_{k-1}}^t \E   | \DD_{k,s}|^4\d s
    +  8K_1^2C_0d^{2}   (1+|x|^4) \eta_k^6,\ \ \ t\in [t_{k-1},t_k].
\end{align*}
By Grownwall's inequality,  we find a constant $C_1\in (0,\infty)$ depending on $K_1,K_2,\eta$ and $\aa$,  such that   \eqref{D1} holds.
  On the other hand, by  $(A1)$ for $\aa\in (1,2]$, we have
\beg{align*}
|\E[\DD_k|Y_{t_{k-1}^x}]|&= |\E[X_{t_{k-1},t_k}(z)- Y_{t_{k-1},t_k}(z)]|_{z=Y_{t_{k-1}^x}} \\
&= \bigg|\int_{t_{k-1}}^{t_k} \E[b\left(X_{t_{k-1},t}(z)\right)-b(z)]\d t\bigg|_{z=Y_{t_{k-1}^x}} \\
&=\bigg|\int_{t_{k-1}}^{t_k} \int_{0}^{1} \E \left[\nn_{ X_{t_{k-1},t}(z)- z}b\left(z+r(X_{t_{k-1},t}(z)-z)\right) \right] \d r\d t\bigg|_{z=Y_{t_{k-1}^x}} \\
&\le \bigg|\int_{t_{k-1}}^{t_k} \int_{0}^{1} \E \left[\nn_{ X_{t_{k-1},t}(z)- z}\left(b\left(z+r(X_{t_{k-1},t}(z)-z)\right)-b(z) \right) \right] \d r\d t\bigg|_{z=Y_{t_{k-1}^x}}\\
&\quad +\bigg|\int_{t_{k-1}}^{t_k} \E[\nn_{ X_{t_{k-1},t}(z)- z}b(z)]\d t\bigg|_{z=Y_{t_{k-1}^x}}  \\
&\le K_1 \int_{t_{k-1}}^{t_k}  \big(\E|X_{t_{k-1},t}(z)- z|^\aa\big)_{z=Y_{t_{k-1}^x}}  \d t
+K_1 \int_{t_{k-1}}^{t_k}  \d t \int_{t_{k-1}}^t \E |b(X_{t_{k-1},s}(z))|_{z=Y_{t_{k-1}^x}} \d s.
\end{align*}
So, by the linear growth of $b$, Lemma \ref{L2.2} and Lemma \ref{C2.3},  we find  constants $k_1,k_2 \in (0,\infty)$ depending on $K_1,K_2,\eta$ and $\aa,$ such that
\begin{align*}
\E |\E[\DD_k|Y_{t_{k-1}^x}]|  & \le  k_1 \eta_k^{1+\ff\aa 2} \E (d^{\frac{\alpha}{2}}+| Y_{t_{k-1}}^x |^\aa ) \\
&\le  k_2 d^{\ff \aa 2} (1+|x|^\aa) \eta_k^{1+\ff\aa 2}.
\end{align*}
Hence,  \eqref{D2} holds.
 By \eqref{D0}, \eqref{D1},  \eqref{D2}    and Lemma \ref{L3.1}, we find constants $p_1,p_2 \in (0,\infty)$ depending only on $K_1,K_2,\eta$  and $\aa$,  such that $\|\nn f\|_\infty\le 1$ implies
  \beg{align*}
  J_1& \le \big\|\|\nn^2 P_{t_k,t}f\|_{\mathrm{op}}\big\|_\infty \E|\DD_k|^2 \\
  & \le p_1d^{\frac{3}{2}} (1+|x|^2) \e^{-\kk_2(t-t_k)} (t-t_k)^{-\ff{1}2} \eta_k^{3},\\
  \end{align*}
  and
 \beg{align*}
J_3&= \big|\E\<\nn P_{t_k, t} f(Y_{t_{k-1}}^x), \E[\DD_k|Y_{t_{k-1}}^x]\>\big| \\&\le \big\|\|\nn  P_{t_k,t}f\|_{\mathrm{op}}\big\|_\infty \E\big|\E[\DD_k|Y_{t_{k-1}}^x]\big|\\
&  \le p_2 d^{\frac{\alpha}{2}}(1+|x|^2) \e^{-\kk_2(t-t_k)}  \eta_k^{1+\ff\aa 2},\ \ t>t_k.\end{align*}
Combining this with \eqref{TT},  we derive the desired estimate provided we find a constant $\kk>0$ such that
\beq\label{FN} J_2   \le  \kk d^\frac{3}{2} (1+|x|^2)\e^{-\kk_2(t-t_k)}  \eta_k^{2 \land \frac{3\alpha}{2}},\ \ t>t_k.  \end{equation}

(c) To verify \eqref{FN}, we apply the integration by parts formula for Malliavin's derivative. It suffices to prove for large $k$, say $k\ge 2$. Let
$$m:=\inf\big\{i\in \Z^+:\ t_k-t_i\le 3\eta_1\big\}.$$ We have $0\le m\le k-2$ and
\beq\label{1} \eta_1\le t_{k-1}-t_m\le 3\eta_1.\end{equation}
Let $\{e_l\}_{1\le l\le d}$ be the canonical orthonormal basis in $\R^d$. For each $1\le l\le d,$ let $h_l(0)=0$,
\beg{align*}
&h_l'(t):= \sum_{i=m}^{k-1} 1_{[t_i,t_{i+1})}(t) \si^{-1}\Big(\ff{e_l}{t_k-t_m} - \ff{t_i-t_m}{t_k-t_m}\nn_{e_l} b(Y_{t_i}^x)\Big),\ \ t\in [0,t_k].\\
\end{align*}
Since $h_{l}(t)$ is adapted with respect to the filtration generated by the Brownian motion, it is clear that (see \cite[Proposition 1.3.11]{N06})
\begin{align*}
    \dd(h_l):= \int_0^{t_k}\<h_l'(t),\d W_2\> = \int_{t_m}^{t_k} \<h_l'(t),\d W_2\>,
\end{align*}
where $\delta$ is the divergence operator in Malliavin calculus. By \eqref{1}, there exist constants $c_1,c_2>0$ depending only on $K_{1},\eta$ such that
\beq
\label{2} |h_l'(t)|\le c_11_{[t_m,t_k]}(t),\ \ \E[\dd(h_l)^4] \le c_1 \E\bigg(\int_{t_m}^{t_k} |h_l'(t)|^2\d t\bigg)^2 \le c_2.
\end{equation}
Let $D_{h_l}$ be the Malliavin derivative along $h_l$. We claim that
\beq\label{2'}D_{h_l}Y_{t_i}^x=\ff{t_i-t_m}{t_k-t_m} e_l,\ \ \ m\le i\le k.\end{equation}
The formula  is trivial for $i=m.$ If it holds for some $m\le i\le k-1,$ by the definition of $Y_t^x$, we have
\beg{align*}
D_{h_l}Y_{t_{i+1}}^x
&=D_{h_l} Y_{t_{i}}^{x}+\eta_{i+1} D_{h_l} b(Y_{t_{i}}^{x})+\sigma\left(h_{l}(t_{i+1})-h_{l}(t_i)\right)\\
&= \ff{t_i-t_m}{t_k-t_m} e_l + \nn_{e_l}b(Y_{t_i}^x) \ff{(t_{i+1}-t_i)(t_i-t_m)}{t_k-t_m} +\si(h_{l}(t_{i+1})-h_{l}(t_i)) \\
&= \ff{t_{i+1}-t_m}{t_k-t_m} e_l,
 \end{align*}
 so that by induction we derive \eqref{2'}, which together with $(A1)$ and \eqref{1}  yields
 \beq\label{2''}   |D_{h_l}(Y_{t}^x-Y_{t_{k-1}}^x)|\le \big|\nn_{D_{h_l}Y_{t_{k-1}}^x}b(Y_{t_{k=1}}^x) (t-t_{k-1}) \big| + |\si (h_l(t)-h_l(t_{k-1}))| \le c_3 \eta_k
 \end{equation}
for some constant $c_3>0$ and all $t\in [t_{k-1},t_k]$.  In particular,
 \beq\label{3} |D_{h_l} \tt \DD_k| =\big|D_{h_l}Y_{t_k}^x - D_{h_l}Y_{t_{k-1}}^x\big|\le c_3 \eta_k.\end{equation}
 Moreover,
 $\DD_{k,t}$  in \eqref{KT} satisfies
 $$D_{h_l}\DD_{k,t}= D_{h_l} X_{t_{k-1},t}(Y_{t_{k-1}}^x) - D_{h_l} Y_{t}^x,$$ so that
 \beg{align*}
 \ff{\d}{\d t} D_{h_l}\DD_{k,t}&= \Big(\nn_{D_{h_l} X_{t_{k-1},t}(Y_{t_{k-1}}^x) }b\Big)\big(X_{t_{k-1},t}(Y_{t_{k-1}}^x)\big) - \Big(\nn_{D_{h_l} Y_{t_{k-1}}^x}b\Big)(Y_{t_{k-1}}^x)\\
 &= \Big(\nn_{D_{h_l} \DD_{k,t}}b\Big)\big(X_{t_{k-1},t}(Y_{t_{k-1}}^x)\big)  + \Big(\nn_{D_{h_l}(Y_{t}^x-Y_{t_{k-1}}^x) }b\Big)\big(X_{t_{k-1},t}(Y_{t_{k-1}}^x)\big)\\
 &\qquad +\Big(\nn_{D_{h_l} Y_{t_{k-1}}^x }b\Big)\big(X_{t_{k-1},t}(Y_{t_{k-1}}^x)\big) - \Big(\nn_{D_{h_l} Y_{t_{k-1}}^x}b\Big)(Y_{t_{k-1}}^x)\\
 &:=\Big(\nn_{D_{h_l} \DD_{k,t}}b\Big)\big(X_{t_{k-1},t}(Y_{t_{k-1}}^x)\big)+\mathcal{R}_{k}(t).
 \end{align*}
 Note that $D_{h_l}\DD_{k,t_{k-1}}=0$. So solving this ODE yields,
 \begin{align*}
      D_{h_l}\DD_{k,t}&=  \int_{t_{k-1}}^{t}\exp\left(\int_{s}^{t} \nabla b\left(X_{t_{k-1},r}(Y_{t_{k-1}}^{x})\right) \d r \right) \mathcal{R}_{k}(s) \mathrm{d}s,\quad \forall t\in [t_{k-1},t_{k}].
 \end{align*}
 Combining this with \eqref{D1}, \eqref{2'}, \eqref{2''} and $(A1)$,  we can find constants $c_4,c_5>0$ such that
  \begin{align}\label{4}
    \begin{split}
        \E|D_{h_l} \DD_k|^2 &= \E|D_{h_l}\DD_{k,t_k}|^2 \\
   & \le \eta_{k}\int_{t_{k-1}}^{t_{k}} \exp\left(2K_{1}\eta_{k} \right) \mathbb{E}\left|\mathcal{R}_{k}(s)\right|^{2} \mathrm{d}s\\
    &\le \eta_{k}\int_{t_{k-1}}^{t_{k}} K_{1}^{2}\exp\left(2K_{1}\eta_{k} \right) \mathbb{E}\left(\left|D_{h_l}(Y_{s}^x-Y_{t_{k-1}}^x)\right|^{2}+\left| \Delta_{k,s}\right|^{2(\alpha-1)}\right) \mathrm{d}s\\
    &\le c_4 \eta_k^2 \Big(\eta_k^2+ \sup_{t\in [t_{k-1},t_k]} \E|\DD_{k,t}|^{2(\aa-1)}\Big)\\
    &\le c_5   d^{\aa-1} (1+|x|^2)^{\aa-1} \eta_k^{2+\min\{2,3(\aa-1)\}} .
    \end{split}
 \end{align}
 Let
 $$\Xi_r:= (1-r)Y_{t_{k-1}}^x+ rY_{t_k}^x,\ \ r\in [0,1].$$
 By \eqref{1} and  \eqref{2'}, it is clear that
 $$D_{h_l} \Xi_r= \Big(\ff{(1-r)(t_{k-1}-t_m)}{t_k-t_m} +r\Big)e_l, $$
 and
 \begin{align*}
    \ff{(1-r)(t_{k-1}-t_m)}{t_k-t_m} +r\ge \ff{ t_{k-1}-t_m}{t_k-t_m}\ge \ff 1 3.
 \end{align*}
 So,  for any function $f$ with $\|\nn f\|_\infty\le 1$, by    \eqref{GR1}  and the integration by parts formula for Malliavin derivative (see, for instance, \cite[Section 1.3]{N06} for more details), we then have
 \begin{align*}
     &\big|\E[\nn_{\tt \DD_k}\nn_{\DD_k} P_{t_k,t}f(\Xi_r)]\big|
     =\left|\sum_{j,l=1}^{d} \mathbb{E}\left[\left(\Delta_{k} \right)_{l} (\tt \Delta_{k} )_{j} \nabla_{e_{l}} \nabla_{e_{j}} P_{t_{k},t}f\left( \Xi_{r}\right) \right]\right|\\
     &\le 3\left|\sum_{j,l=1}^{d} \mathbb{E} \left[ \left(\Delta_{k} \right)_{l} (\tt \Delta_{k} )_{j} \nabla_{D_{h_l}\Xi_r} \nabla_{e_{j}} P_{t_{k},t}f\left( \Xi_{r}\right) \right]\right|\\
     &=3\left|\sum_{j,l=1}^{d} \mathbb{E} \left[ \left(\Delta_{k} \right)_{l} (\tt \Delta_{k} )_{j} D_{h_l}\left\{\nabla_{e_{j}} P_{t_{k},t}f\left( \Xi_{r}\right) \right\} \right]\right|\\
     &=3\left|\sum_{j,l=1}^{d} \mathbb{E} \left[\delta \left(\left(\Delta_{k} \right)_{l} (\tt \Delta_{k} )_{j} h_{l}\right)\nabla_{e_{j}} P_{t_{k},t}f\left( \Xi_{r}\right)  \right]\right|\\
     &=3\left|\sum_{j,l=1}^{d} \mathbb{E} \left[\left(-D_{h_{l}}\{\left(\Delta_{k} \right)_{l}\} (\tt \Delta_{k} )_{j} -\left(\Delta_{k} \right)_{l} D_{h_{l}}\{(\tt \Delta_{k} )_{j}\}+\left(\Delta_{k} \right)_{l} (\tt \Delta_{k} )_{j} \delta(h_{l})  \right)\nabla_{e_{j}} P_{t_{k},t}f\left( \Xi_{r}\right)  \right]\right|\\
     &=3\left| \mathbb{E} \left[\nabla_{\tt \Delta_{k}} P_{t_{k},t}f\left( \Xi_{r}\right) \left(\sum_{l=1}^{d}D_{h_{l}}\{\left(\Delta_{k} \right)_{l}\} \right)\right]\right|+3\left| \mathbb{E} \left[ \sum_{l=1}^{d}\left(\Delta_{k} \right)_{l}\left(\nabla_{D_{h_{l}}\{ \tt \Delta_{k}\}} P_{t_{k},t}f\left( \Xi_{r}\right)\right)\right]\right|\\
     &\quad +3\left|\mathbb{E}\left[\nabla_{\tt \Delta_{k}} P_{t_{k},t}f\left( \Xi_{r}\right) \left(\sum_{l=1}^{d}  \left(\Delta_{k} \right)_{l} \delta(h_{l}) \right)\right] \right|\\
     &:=T_{1}+T_{2}+T_{3},
 \end{align*}
 where $(\Delta_{k})_{j}$ denotes the $j$-th component of $\Delta_{k}\in \mathbb{R}^{d}$. It follows from Cauchy-Schwarz inequality, \eqref{GR1}, \eqref{D0}, \eqref{2}, \eqref{3} and \eqref{4} that
 \begin{align*}
     T_{1}&\le 3\kk_{1}\e^{-\kk_2 (t-t_k)} \left(\mathbb{E}|\tt \Delta_{k}|^{2}\right)^{\frac{1}{2}} \left(\mathbb{E}\left(\sum_{l=1}^{d}D_{h_{l}}\{\left(\Delta_{k} \right)_{l}\} \right)^{2} \right)^{\frac{1}{2}}\\
     &\le c_{6}\e^{-\kk_2 (t-t_k)} d^{\frac{\alpha+1}{2}} \left(1+|x|^{2}\right)^{\frac{\alpha}{2}} \eta_{k}^{\frac{5}{2}\land \frac{3\alpha}{2}},
 \end{align*}
 \begin{align*}
     T_{2}&\le 3\kk_{1}\e^{-\kk_2 (t-t_k)}\left(\mathbb{E}| \Delta_{k}|^{2}\right)^{\frac{1}{2}} \left(\mathbb{E}\sum_{l=1}^{d}\left(D_{h_{l}}\{\tt \Delta_{k} \} \right)^{2} \right)^{\frac{1}{2}} \\
     & \le c_{6}\e^{-\kk_2 (t-t_k)} d(1+|x|^{2})^{\frac{1}{2}} \eta_{k}^{\frac{5}{2}},
 \end{align*}
 and
 \begin{align*}
     T_{3}&\le 3\kk_{1}\e^{-\kk_2 (t-t_k)} \left(\mathbb{E}|\tt \Delta_{k}|^{2}\right)^{\frac{1}{2}} \left(\mathbb{E} |\Delta_{k}|^{4} \right)^{\frac{1}{4}} \left(d\sum_{l=1}^{d}\mathbb{E}\left[ \delta(h_{l})^{4}\right] \right)^{\frac{1}{4}}\\
     &\le c_{6}\e^{-\kk_2 (t-t_k)} d^{\frac{3}{2}} (1+|x|^{2})\eta_{k}^{2},
 \end{align*}
 for some constant $c_{6}>0$. Combining these estimates, we derive \eqref{FN}.
\end{proof}

\beg{proof}[Proof of Theorem \ref{T2}] For $p\in (0,1)$, according to \cite[Lemma 2.1]{HRW} for $\psi(r)=r^p$ and $\W_0(\mu_1,\mu_2)=\ff 1 2\|\mu_1-\mu_2\|_{\mathrm{TV}}$, it holds for any two probability measures $\mu$ and $\nu$ that
\begin{align*} 
    \W_{p}(\mu,\nu)&\le \inf_{t>0} \big\{2\ss d t^{\ff p 2} \W_0(\mu,\nu) + d t^{\ff{p-1}2} \W_1(\mu,\nu)\big\}\\
    &
        \le 2^{2-p} d^{\ff{1+p}2} \W_0(\mu,\nu)^{1-p} \W_1(\mu,\nu)^{p}.
\end{align*}
So, we only need to prove statements (1) and (2).

It is well known that $P_t$ is ergodic and admits a unique invariant measure $\mu$ under $(A1)$ (see for instance \cite{Ha10}). By \eqref{GR1} and \eqref{GR2}, we have
\begin{align}\label{ergodic}
        \W_1(\L(X_t^x),\mu)
    &=\sup_{\|\nabla f\|_{\infty}<1} \left|P_{t}f(x)-\int_{\R^{d}} P_{t}f(y) \mu(\mathrm{d}y) \right| \notag\\
    &\le \sup_{\|\nabla f\|_{\infty}<1} \left|\int_{\R^{d}} \left(P_{t}f(x)-P_{t}f(y)\right) \mu(\mathrm{d}y) \right| \notag\\
    &\le \sup_{\|\nabla f\|_{\infty}<1} \|\nabla P_{t}f\|_{\infty}\int_{\R^{d}} |x-y| \mu(\mathrm{d}y) \notag\\
    &\le \kk_1 \e^{-\kk_2 t} \mu(|x-\cdot|)\le c_1 (d^{\frac{1}{2}}+|x|) \e^{-\kk_2 t},
 \end{align}
{where the last inequality is a consequence of the ergodicity and Lemma \ref{L2.2} } and similarly,
\begin{align*}
    \W_0(\L(X_t^x),\mu)
    &=\sup_{\|f\|_{\infty}<1} \left|P_{t}f(x)-\int_{\R^{d}} P_{t}f(y) \mu(\mathrm{d}y) \right|\\
    &\le \min\big\{1, \kk_1 \e^{-\kk_2 t}t^{-\ff 1 2} \mu(|x-\cdot|)\big\} \\
    & \le c_1 (d^{\frac{1}{2}}+|x|) \e^{-\kk_2 t},\ \ \ \ \ \ \ \ t>0,
\end{align*}
for some constant $c_1=c_1(K_1,K_2)\in (0,\infty).$
So we only need to prove the desired upper bound for $\W_i(\L(X_{t_n}^x), \L(Y_{t_n}^x)), i=1,0.$ For any $z\in \R^d$, let $\mu_n^z:=\L(X_{t_{n-1},t_n}(z)), \nu_n^z:= \L(X_{t_{n-1},t_n}(z)).$

(a) Estimate  on $\W_1$. When $\aa\in (0,1]$, \eqref{W1} follows from \eqref{SM} and Lemma \ref{L3.2}(1). When $\aa\in (1,2]$, by \eqref{SM} and Lemma \ref{L3.2}(2) for $t=t_n$,
we find a constant $\kappa=\kappa(K_1,K_2,\eta,\aa)\in (0,\infty)$ such that for any $1\le k\le n-1,$
\beg{align*}
\W_1(\L(X_{t_n}^x,\L(Y_{t_n}^x) )  &\le  \kk d^{\frac{3}{2}} (1+|x|^2)  \sum_{k=1}^{n-1} \e^{-\kappa_{2} (t_n-t_k)}  \eta_k^{1+\ff\aa 2} \\
&\quad + \sup_{\|\nn f\|_\infty\le 1} |Q_{0,t_{n-1}}(P_{t_{n-1},t_n}-Q_{t_{n-1},t_n})f(x)|.\end{align*}
Notice that
\begin{align} \label{inequ}
    \sum_{k=1}^{n-1} \e^{-\kappa_{2} (t_n-t_k)}  \eta_k^{1+\ff\aa 2}
    \ge \eta_{n}^{\frac{\alpha}{2}} \sum_{k=1}^{n-1} \e^{-\kappa_{2} (t_n-t_k)}  \eta_k
    \ge \eta_{n}^{\frac{\alpha}{2}} \int_{0}^{t_{n-1}} e^{-\kappa_{2}(t_{n}-t)} \mathrm{d}t\ge \frac{1}{\kappa_{2}}\eta_{n}^{\frac{\alpha}{2}}.
\end{align}
At the same time, (\ref{D1}) implies that, for any $f$ satisfying $\|\nabla f\|_{\infty}\le 1$,
\begin{align*}
    \left|Q_{0,t_{n-1}}\left(P_{t_{n-1},t_{n}}-Q_{t_{n-1},t_{n}}\right)f(x)\right|
    &=\left|\mathbb{E}\left[f\left(X_{t_{n-1},t_{n}}(Y_{t_{n-1}}^{x})\right)-f(Y_{t_{n}}^{x}) \right]  \right|\\
    &\le \mathbb{E}\left| \Delta_{n}\right| \le C_{1} d^{\frac{1}{2}}(1+|x|^{2})^{\frac{1}{2}} \eta_{n}^{\frac{3}{2}},
\end{align*}
where $\Delta_{n}$ is defined as in (\ref{D}). Thus, according to (\ref{inequ}), the last term can be absorbed by the sum term and \eqref{W1'} holds for a possibly larger constant $\kk \in (0,\infty)$.\par
(b) Estimate on $\W_0$. By  \eqref{GR2} and Lemma \ref{L3.2} for $t=t_n$, we can find  a constant $k_0=k_0(K_1,K_2,\eta,\aa)\in (0,\infty)$ such that for any $1\le k\le n-1,$
\beg{align*}
&\ \ \ \ \ \sup_{\|f\|_\infty\le 1} |Q_{0,t_{k-1}}(P_{t_{k-1}, t_k}- Q_{t_{k-1},t_k})P_{t_k, t_n}f(x)|\\
&=\sup_{\|f\|_\infty\le 1} |Q_{0,t_{k-1}}(P_{t_{k-1}, t_k}- Q_{t_{k-1},t_k})P_{t_k, \frac{t_k+t_n}{2}}\{P_{\frac{t_k+t_n}{2}, t_n}f\}(x)|\\
&\le   \left(\sup_{\|f\|_{\infty}\le 1}\|\nabla P_{\frac{t_k+t_n}{2}, t_n}f\|_{\infty} \right)\sup_{\|\nn g\|_\infty\le 1} |Q_{0,t_{k-1}}(P_{t_{k-1}, t_k}- Q_{t_{k-1},t_k}) P_{t_k,\frac{(t_k+t_n)}{2}}g(x)|\\
& \le 2\kk_1 \e^{-\kk_2(t_n-t_k)/2}(t_n-t_k)^{-\ff 1 2} \sup_{\|\nn g\|_\infty\le 1} |Q_{0,t_{k-1}}(P_{t_{k-1}, t_k}- Q_{t_{k-1},t_k}) P_{t_k,\frac{(t_k+t_n)}{2}}g(x)|\\
&\le   1_{(0,1]}(\aa) \ff{k_0 d^{\ff 1 2} (1+|x|)}{\ss{t_n-t_k}} \e^{-\kk_2(t_n-t_k)}   \eta_k^{1+\ff\aa 2}+ 1_{(1,2]}(\aa)\ff{k_0 d^{\frac{3}{2}} (1+|x|^2)}{\ss{t_n-t_k}}  \e^{-\kk_2(t_n-t_k)} \eta_k^{1+\ff\aa 2}.
\end{align*}
Combining this with  (\ref{SM}), we derive \eqref{W0} and \eqref{W0'}  provided  there exists a positive constant $k_{1}=k_{1}(K_1,K_2,\eta,\aa)$ such that
 \beq\label{LST} \sup_{\|f\|_\infty\le 1} |Q_{0,t_{n-1}}(P_{t_{n-1}, t_n}- Q_{t_{n-1},t_n}) f(x)|\le k_{1} (1+|x|)d^{\ff 1 2}\eta_n^{\ff {1  +(\aa\land 1)}2},\ \ n\ge 2.
 \end{equation}
To this end, for any $z\in\R^d$, let
$$\mu_n^z=\L(X_{t_{n-1},t_n}(z)), \ \ \ \nu_n^z=\L(Y_{t_{n-1},t_n}(z)).$$ By Lemma \ref{L2.2}, \eqref{LST} follows if  we can   find constants $k_2 =k_2 (K_1,K_2,\eta,\aa) \in (0,\infty)$   such that {
\beq\label{LST2}  \W_0(\mu_n^z,\nu_n^z)\le  k_2 (d^{\ff 1 2}+|z|) \eta_n^{ \frac{1+(\aa\land 1)}{2} },\ \ n\ge 2, z\in\R^d.\end{equation}}

Let us now show (\ref{LST}). Write
$$\d X_{t_{n-1},t}(z) = b(z)\d t+ \si \d\tt B_t,\ \ \ X_{t_{n-1},t_{n-1}}(z)=z,\ t\in [t_{n-1},t_n],$$
where
$$\tt B_t= B_t -\int_{t_{n-1}}^{t} \big\{b(z)- b(X_{t_{n-1},s}(z) )\big\}\d s,\ \ \ t\in [t_{n-1},t_n].$$
Let
$R:= \exp\left(\int_{t_{n-1}}^{t_{n}}\<b(z)- b\left(X_{t_{n-1},s}(z) \right),\d B_s\> -\ff 1 2\int_{t_{n-1}}^{t_{n}} \big|b(z)- b\left(X_{t_{n-1},s}(z)\right)\big|^2\d s\right).$
In order to apply the Girsanov's theorem, we first show that $\mathbb{E}[R]=1$. By the continuity assumption of $b$ in $(A1)$ and Young's inequality, there exist a positive constant $k_{3}=k_{3}(K_{1},\alpha)$ such that
\begin{align} \label{bb}
    \left| b(X_{t_{n-1},s}(z) )-b(z)\right|^{2}\le 4K_{1}^{2}\left|X_{t_{n-1},s}(z)-z \right|^{2\land 2\alpha}\le k_{3}\left(\left|X_{t_{n-1},s}(z)-z \right|^{2}+1\right).
\end{align}
So, according to \cite[Ch.8,Exercise 1.40]{RY13}(another Novikov's type criterion), it suffices to show that,
\begin{align}\label{Novikov}
    \mathbb{E}\left[\exp\left(a \left|X_{t_{n-1},t}(z) -z \right|^{2} \right) \right]\le c,
\end{align}
for any $t\in [t_{n-1},t_{n}]$ and two constants $a$ and $c$. It follows from Itô's formula, $(A1)$ and Young's inequality that, for any $t\in [t_{n-1},t_{n}]$,
\begin{align*}
    \left|X_{t_{n-1},t}(z)-z \right|^{2}
    &=2\int_{t_{n-1}}^{t} \langle X_{t_{n-1},s}(z)-z,b(X_{t_{n-1},s}(z) ) \rangle \mathrm{d}s+2\int_{t_{n-1}}^{t} \langle X_{t_{n-1},s}(z)-z, \sigma \mathrm{d}B_{s} \rangle \\
     &\quad +\|\sigma\|_{\mathrm{HS}}^{2}(t-t_{n-1})\\
    &\le -2K_{2}\int_{t_{n-1}}^{t} \left|X_{t_{n-1},s}(z)-z\right|^{2}  \mathrm{d}s+2\int_{t_{n-1}}^{t} \langle X_{t_{n-1},s}(z)-z, \sigma \mathrm{d}B_{s} \rangle\\
    &\quad + 2\int_{t_{n-1}}^{t} \left|X_{t_{n-1},s}(z)-z\right| \left| b(z)\right|\mathrm{d}s+K_{1}d\eta_{n}\\
    &\le -2\Tilde{K}_{2}\int_{t_{n-1}}^{t} \left|X_{t_{n-1},s}(z)-z\right|^{2}  \mathrm{d}s+2\int_{t_{n-1}}^{t} \langle X_{t_{n-1},s}(z)-z, \sigma \mathrm{d}B_{s} \rangle+k_{4}d\eta_{n}
\end{align*}
for some positive constant $k_{4}=k_{4}(K_{1},K_{2},b(z))$ and $\Tilde{K}_{2}\in (0,K_{2})$. As a consequence, for any $\gamma>0$, we have
\begin{align*}
    &\mathbb{E}\left[\exp\left(\gamma \left|X_{t_{n-1},t}(z)-z \right|^{2}+2\Tilde{K}_{2}\gamma\int_{t_{n-1}}^{t} \left|X_{t_{n-1},s}(z)-z\right|^{2}  \mathrm{d}s \right) \right]\\
    &\le e^{k_{4}d\eta_{n}}\mathbb{E}\left[\exp\left( 2\gamma\int_{t_{n-1}}^{t} \langle X_{t_{n-1},s}(z)-z, \sigma \mathrm{d}B_{s} \rangle\right)  \right].
\end{align*}
Further more, Hölder’s inequality , $(A1)$ and the local exponential martingale property implies that, for any fixed $\gamma\in (0,\frac{\Tilde{K}_{2}}{4K_{1}^{2}})$ and $t\in [t_{n-1},t_{n}]$,
\begin{align*}
    &\mathbb{E}\left[\exp\left( 2\gamma\int_{t_{n-1}}^{t} \langle X_{t_{n-1},s}(z)-z, \sigma \mathrm{d}B_{s} \rangle\right)  \right]\\
    &\le \left(\mathbb{E}\left[\exp\left( 4\gamma\int_{t_{n-1}}^{t} \langle X_{t_{n-1},s}(z)-z, \sigma \mathrm{d}B_{s} \rangle-8\gamma^{2}\int_{t_{n-1}}^{t}  \left|\sigma\left(X_{t_{n-1},s}(z)-z\right) \right|^{2}  \mathrm{d}s \right)  \right] \right)^{\frac{1}{2}}\\
    &\quad \times \left( \mathbb{E}\left[\exp\left( 8\gamma^{2}\int_{t_{n-1}}^{t}  \left|\sigma\left(X_{t_{n-1},s}(z)-z\right) \right|^{2}  \mathrm{d}s \right)  \right] \right)^{\frac{1}{2}}\\
    &\le \left( \mathbb{E}\left[\exp\left( 8\gamma^{2}\int_{t_{n-1}}^{t}  \left|\sigma\left(X_{t_{n-1},s}(z)-z\right) \right|^{2}  \mathrm{d}s \right)  \right] \right)^{\frac{1}{2}}\\
    &\le \left( \mathbb{E}\left[\exp\left( 8K_{1}^{2}\gamma^{2}\int_{t_{n-1}}^{t}  \left|X_{t_{n-1},s}(z)-z \right|^{2} \mathrm{d}s \right)  \right] \right)^{\frac{1}{2}}\\
    &\le \left( \mathbb{E}\left[\exp\left(\gamma \left|X_{t_{n-1},t}(z)-z \right|^{2}+ 2\Tilde{K}_{2}\gamma\int_{t_{n-1}}^{t}  \left|X_{t_{n-1},s}(z)-z \right|^{2}  \mathrm{d}s \right)  \right] \right)^{\frac{1}{2}}.
\end{align*}
Hence, (\ref{Novikov}) follows with $a=\gamma$, $c=e^{2k_{4}d\eta_{n}}$ and $R$ is a martingale.\par
Now, by Girsanov's theorem, $(\tt B_t)_{t\in [t_{n-1},t_n]} $ is a Brownian motion under the probability measure $R\d\P$.
Clearly, $\mathbb{E}[R-1]=0$. {Notice that, by Cauchy-Schwarz inequality,
\begin{align*}
    \left(\mathbb{E}|R-1|\right)^{2}
    =\left(\mathbb{E}\left[\frac{|R-1|}{\sqrt{1+\frac{R-1}{3}}} \sqrt{1+\frac{R-1}{3}} \right]\right)^{2}
    \le  \mathbb{E}\left[ \frac{(R-1)^{2}}{1+\frac{R-1}{3}} \right] \mathbb{E}\left[1+\frac{R-1}{3} \right].
\end{align*}
}
Combining this with the fact that
$$(1+x)\log(1+x)-x\ge \frac{1}{2}\left(\frac{x^{2}}{1+\frac{x}{3}}\right),\, \forall x> -1,$$
yields
\begin{align*}
    \frac{1}{2}\left(\mathbb{E} |R-1|\right)^{2}
    =\frac{1}{2} \frac{\left(\mathbb{E}|R-1|\right)^{2}}{\mathbb{E}\left[1+\frac{R-1}{3}\right]}
    \le \frac{1}{2}\mathbb{E}\left[ \frac{(R-1)^{2}}{1+\frac{R-1}{3}} \right]
    \le \mathbb{E}[R\log R],
\end{align*}
which is also known as a Pinsker type inequality.
Hence, combining this with the definition of $\W_0$, Girsanov's theorem   and Lemma \ref{C2.3}, we find $k_5=k_5(K_1,K_2,\eta,\aa)\in (0,\infty)$ such that
\beg{align*} &\W_0(\mu_n^z,\nu_n^z)^2=\ff 1 4 \sup_{\|f\|_\infty\le 1} |\mu_n^z(f)-\nu_n^z(f)|^2 \\
&= \ff 1 4 \sup_{\|f\|_\infty\le 1} \big| \E[f(X_{t_{n-1},t_n}(z))- R f(X_{t_{n-1},t_n}(z))]\big|^2
  \le \ff 1 4 \left[\E|R-1|\right]^2 \le  \frac{1}{2}\E [R\log R] \\
  & = \ff 1 4 \int_{t_{n-1}}^{t_n} \E[R |b(z)- b\left(X_{t_{n-1},s}(z)\right)|^2]\d s
 = \ff 1 4 \int_{t_{n-1}}^{t_n}\E[  |b(z)- b\left(Y_{t_{n-1},s}(z) \right)|^2]\d s  \\
 & \le  K_1^2  \int_{t_{n-1}}^{t_n} \E\big[|Y_{t_{n-1},s}(z)-z|^2+ |Y_{t_{n-1},s}(z)-z|^{2(1\land\aa)}\big]\d s
 \le c_2 (d+|z|^2) \eta_n^{1+ (1\land\aa)}.\end{align*}
So,  \eqref{LST2} holds.

    (c) For the particular case, let $\eta_{k}=\frac{\theta}{k}$ for some $\theta>\frac{\alpha}{2c_{2}}$. Since $\theta \left(\log(n)-\log(k+1)\right)\le t_{n}-t_{k}\le \theta\left( \log(n-1)-\log(k) \right)$, we can find a positive constant $k_{6}=k_{6}(\theta, \alpha,c_{2})$ (may vary from line to line) such that
    \begin{align*}
        \sum_{k=[\frac{n-1}{2}]+1}^{n-1}    \e^{ -c_{2}(t_{n}-t_{k}) } (t_n-t_k)^{-\frac{1}{2}} \eta_{k}^{1+\ff \aa 2}
        &\le k_{4}   \sum_{k=[\frac{n-1}{2}]+1}^{n-1} \left( \frac{k}{n}\right)^{\theta c_{2}} \frac{1}{\sqrt{\log n-\log k}} k^{-1-\frac{\alpha}{2}}\\
        &\le k_{4} n^{-\frac{\alpha}{2}} \sum_{k=1}^{n-1}\frac{1}{\sqrt{\log n-\log(k)}} k^{-1}\\
        &\le k_{4} n^{-\frac{\alpha}{2}} \int_{\frac{n-1}{2}}^{n-1} \frac{1}{t\sqrt{\log n-\log(t)}} \mathrm{d}t\\
        &\le k_{4}n^{-\frac{\alpha}{2}},
    \end{align*}
  and
  \begin{align*}
      \sum_{k=1}^{[\frac{n-1}{2}]}\e^{ -c_{2}(t_{n}-t_{k}) } (t_n-t_k)^{-\frac{1}{2}} \eta_{k}^{1+\ff \aa 2}
      &\le k_{4} n^{-\theta c_{2}} \sum_{k=1}^{[\frac{n-1}{2}]} k^{-1+\theta c_{2}-\frac{\alpha}{2}}
      \le k_{4}n^{-\frac{\alpha}{2}}.
  \end{align*}
  So  (\ref{W'}) follows by combining above estimates with (\ref{ergodic}). (\ref{W}) can be proved through the same argument and the proof is complete.
\end{proof}

\section{$\W_p$-estimate for $p>1$: the uniform dissipation case}

To cover typical time dependent models, see Example 3.1 below, we consider the following time in-homogenous   SDE:
 \begin{align} \label{SDE'}
    \d X_t=b_t ( X_t) \d t+\sigma_t  (X_t)\d B_t, \ \ t\ge 0,
\end{align}
where $b:\mathbb{R}^d \times [0,\infty)\mapsto \mathbb{R}^d,\ \ \sigma: \R^d \times [0,\infty) \to\mathbb{R}^{d \times d}$ are measurable. 
The associated continuous time Euler Scheme is defined by
\begin{align} \label{EM'}
    Y_{t}=Y_{t_{k}}+(t-t_{k})b_{t_k}(Y_{t_{k}})+\sigma_{t_{k}}(Y_{t_{k}})\left( B_{t}-B_{t_{k}}\right), \qquad t\in [t_{k},t_{k+1}), k\ge 0, Y_0=X_0.
\end{align}  \par

\subsection{Main result and an example of bridge regression}

\beg{enumerate} \item[$(A2)$] Let $\alpha\in (0,1]$ and $p\in (1,\infty)$. There exist positive constants $K_{1},K_{2}$   such that for any $x,y\in \mathbb{R}^{d},s,t\in [0,\infty)$
    \begin{equation} \label{b}
        \left|b_t(x)-b_s(y) \right|\le K_{1}\left( \left| x-y\right|^{\alpha}+|x-y| +|t-s|\right),
    \end{equation}
    \begin{equation} \label{sigma}
        \|\sigma_{t}(x)-\sigma_{s}(y)\|_{\mathrm{HS}}\le K_{1}\left(|x-y|+|t-s|\right),\ \ \|\si_{t}(x)\|_{\mathrm{op}}\le K_1,
    \end{equation}
    \begin{align} \label{diss}
    &p\langle b_t(x)-b_t(y),x-y \rangle+\ff p 2 \|\sigma_{t}(x)-\sigma_{t}(y)\|^{2}_{\mathrm{HS}}\\
    &\quad + \ff{p(p-2)}{2|x-y|^2} \big|(\sigma_{t}(x)-\sigma_{t}(y))^*(x-y)\big|^2 \le -K_2|x-y|^{2}.\notag
    \end{align}
    \end{enumerate}
We call (\ref{diss}) the uniform dissipation condition.
\beg{thm}\label{T1}
Assume $(A0)$ and $(A2)$. Then $\eqref{SDE'}$ is well-posed . Moreover, for any $K_2'\in (0,K_2)$, there exists a positive constant $c=c(K_1,K_2, \eta,K_2',p,\aa)$ such that
\beq\label{MU0}
\beg{split}&\W_p(\L(X_{t_n}^x), \L(Y_{t_n}^x))^p
\le  cd^{\ff p 2} (1+|x|^p) \sum_{k=1}^n\eta_k^{1+\ff {\aa p} 2 } \e^{-K_2'(t_n-t_k)},\ \ n\ge 1,\forall x\in\R^d.\end{split}
\end{equation}
Consequently, when $\eta_k=\ff \theta k$ for some constant $\theta\in (0,\infty),$ there exists a positive constant $c'=c'(K_1,K_2,\theta, K_2',p,\aa)$ such that
\beq\label{MU}  \beg{split}&\W_p(\L(X_{t_n}^x), \L(Y_{t_n}^x))^p
\le   c'd^{\ff p 2}  (1+|x|^p) n^{-\big((\theta K_2') \land \frac{\alpha p}{2}\big)},\ \ n\ge 1, \forall x\in \R^d.\end{split} \end{equation}
In particular, if $(\sigma_{t},b_t) =(\si,b)$ does not depend on $t$, then the solution of (\ref{SDE'}) is exponentially ergodic with  unique invariant probability measure $\mu$, and 
\begin{align*}
    &\W_p(\L(X_{t_n}^x), \L(Y_{t_n}^x))^p+\W_p(\mu, \L(Y_{t_n}^x))^p\\
&\le  cd^{\ff p 2} (1+|x|^p) \sum_{k=1}^n\eta_k^{1+\ff {\aa p} 2 } \e^{-K_2'(t_n-t_k)},\ \ n\ge 1,\forall x\in\R^d.
\end{align*}
\end{thm}

As an application of Theorem \ref{T1}, we  consider the following optimization problem that arises in the Bridge regression with the shrinkage parameter $\gamma\in (1,2]$ and the tuning parameter $\lambda\ge 0$ (see, i.e., \cite{Fu98} for more details):
\begin{align}\label{tbeta}
    \Tilde{\beta}=\operatorname{argmin}_{\beta\in \mathbb{R}^{d}}\left\{ L(\beta)\right\},L(\beta):=\sum_{i=1}^{N}\left(y_{i}-x_{i}^{T}\beta \right)^{2}+\lambda \sum_{j=1}^{d}\left|\beta_{j}\right|^{\gamma}
\end{align}
where $(x_{i},y_{i})_{\{1\le i\le N\}}$ are $\mathbb{R}^{d}\times \mathbb{R}$-valued data points. In particular, when $\gamma=1,2$, $\Tilde{\beta}$ corresponds to the estimators of the well-known Lasso and Ridge regression, respectively. In practical applications, such optimization problems are usually solved by gradient descent algorithm and its variants. One common variant is the gradient descent algorithm with slowly decreasing Gaussian noise given by the following iterative formula
\begin{align} \label{betagd}
        \beta_{k+1}=\beta_{k}-\eta_{k+1}\nabla L(\beta_{k})+\sqrt{\eta_{k+1}} \sigma_{k+1}\zeta_{k+1},\quad k\ge 0,
\end{align}where $\si_k\downarrow 0$ as $k\uparrow\infty$, and $\{\eta_k\}_{k\ge 1}$ are i.i.d. $d$-dimensional standard Gaussian random variables. Similarly, we may consider it as an approximation of the following SDE:
\begin{align} \label{betaSDE}
    \mathrm{d}\Bar{\beta}_{t}=-\nabla L(\Bar{\beta}_{t})\mathrm{d}t+\Bar{\sigma}_{t}\mathrm{d}B_{t},\quad \Bar{\beta}_{0}=\beta_{0}.
\end{align}
 It has been shown in \cite{G-M93} that, under appropriate assumptions on the drift coefficient, if we set $\sigma_{k}=\frac{1}{\sqrt{\log\log k}} $ for large $k$ and $\eta_{k}=\frac{\theta}{k}$ (correspondingly, $\Bar{\sigma}_{t}\sim \frac{1}{\sqrt{\log t}}$ for large $t$), then this algorithm will converge in probability to $\Tilde{\beta}$. 

It is easy to verify that $L$ is strongly convex, i.e.,
    \begin{align}\label{dissL}
        \langle \nabla L(\beta_{1})-\nabla L(\beta_{2}),\beta_{1}-\beta_{2} \rangle\ge  K |\beta_{1}-\beta_{2}|^{2}, \quad \forall \beta_{1},\beta_{2}\in \R^{d} 
    \end{align}
holds for some constant  $K =K (x,\gamma,\lambda)>0.$  Below we use Theorem \ref{T1} to analyze the convergence rate of this algorithm.  

\textbf{Example 3.1.} 
    Let $\Tilde{\beta},\beta_{k}$ and $\Bar{\beta}_{t}$ be defined as in (\ref{tbeta}), (\ref{betagd}) and (\ref{betaSDE}) respectively.
    For any $p\ge 1$ and $K'\in (0,K)$, let  $\sigma_{k}=\{k^{-\frac{K'\theta}{p}}(\theta \log k)^{-\frac{2}{p}}\wedge 1\}$ and
 $\eta_{k}=\frac{\theta}{k}$ for some $\theta>\frac{(\gamma-1)p}{2K'}$ with $\gamma>1$. Then there exists $c>0$ such that for large n,
     \begin{align*}
        \E|\beta_{n} -\tt\beta |^{p}\le cd^{\frac{p}{2}}\left[|\beta_{0}-\Tilde{\beta}|^p+(1+|\beta_{0}|^p)\right]n^{-\frac{(\gamma-1)p}{2}}.
    \end{align*}

\begin{proof}   By Jensen's inequality, it suffice to consider $p\ge 2$. To apply Theorem \ref{T1}, let $ \bar \sigma_{t}=\{\e^{-\frac{K'}{p}t}t^{-\frac{2}{p}}\wedge 1\}$.
  Consider the gradient flow $\Tilde{\beta}_{t}$ defined by
    \begin{align}\label{gradflow}
        \mathrm{d}\Tilde{\beta}_{t}=-\nabla L(\Tilde{\beta})\mathrm{d}t,\quad \Tilde{\beta}_{0}= {\beta_0}.
    \end{align}
    By It\^o's formula, (\ref{dissL}) and Young's inequality, for any $K'\in (0,K)$, we can find a positive constant $c$ such that
    \begin{align*}
        &\mathrm{d}|\Bar{\beta}_{t}-\Tilde{\beta}_{t}|^{p}-\mathrm{d}M_{t}\\
        &=\left[p|\Bar{\beta}_{t}-\Tilde{\beta}_{t}|^{p-2}\langle \Bar{\beta}_{t}-\Tilde{\beta}_{t},-\nabla L(\Bar{\beta}_{t})+\nabla L(\Tilde{\beta}_{t}) \rangle +\frac{1}{2}p(p-2+d)\bar\sigma_{t}^{2}|\Bar{\beta}_{t}-\Tilde{\beta}_{t}|^{p-2}\right]\mathrm{d}t\\
        &\le \left[-K'|\Bar{\beta}_{t}-\Tilde{\beta}_{t}|^{p}+cd^{\frac{p}{2}}\bar\sigma_{t}^{p}\right]\mathrm{d}t,
    \end{align*}
    which yields that
    \begin{align*}
        \mathbb{E}|\Bar{\beta}_{t}-\Tilde{\beta}_{t}|^{p}
        &\le e^{-K't}\left[|\Bar{\beta}_{0}-\Tilde{\beta}|^p+cd^{\frac{p}{2}}\int_{0}^{t} e^{K_{2}'s} \sigma_{s}^{p}  \mathrm{d}s\right]\\
        &\le  e^{-K't}\left[|\Bar{\beta}_{0}-\Tilde{\beta}|^p+cd^{\frac{p}{2}}\left(\int_{0}^{1} e^{K_{2}'s}   \mathrm{d}s+\int_{1}^{t} \frac{1}{s^{2}}   \mathrm{d}s\right)\right]\\
        &\le e^{-K't}\left[|\Bar{\beta}_{0}-\Tilde{\beta}|^p+c'd^{\frac{p}{2}}\right].
    \end{align*}
    As $t_{n}\sim \theta\log(n)$ for large $n$ and $K'\theta>\frac{(\gamma-1)p}{2}$, it follows that
    \begin{align*}
        \W_{p}(\scr{L}(\Bar{\beta}_{t_{n}}),\delta_{\Tilde{\beta}})^{p}= \mathbb{E}|\Bar{\beta}_{t_{n}}-\Tilde{\beta} |^{p}
        \le \left[|\Bar{\beta}_{0}-\Tilde{\beta}|^p+c'd^{\frac{p}{2}}\right]n^{-K_{2}'\theta}\le \left[|\Bar{\beta}_{0}-\Tilde{\beta}|^p+c'd^{\frac{p}{2}}\right]n^{-\frac{(\gamma-1)p}{2}}.
    \end{align*}
    On the other hand, It can be verified that $L$ and $\bar\sigma_{t}$ satisfy $(A2)$ with $\alpha=\gamma-1$. Hence, Theorem \ref{T1} implies that, for any $K'\in (0,K), n\ge 1$,
    \begin{align*}
    \W_p(\L(\Bar{\beta}_{t_{n}}), \L(\beta_{n}))^p
     \le   c d^{\ff p 2}  (1+|\beta_{0}|^p) n^{- \frac{(\gamma-1) p}{2}}.
    \end{align*}
    Using triangle inequality to combining above two upper bounds gives us the desired result.
\end{proof}

 
\subsection{Proof of Theorem \ref{T1} } 
\ \newline
{Similar to the previous section, we first present the following two lemmas regarding the moment estimates which can be proved through the same way as Lemma \ref{L2.2} and \ref{C2.3}.  }
\beg{lem}\label{L2.2(2)}  Assume that $(A0)$ holds, $\eqref{SDE'}$ is well-posed and there exist positive constants $\kk_1,\kk_2$ such that
\beq\label{ZZ(2)} \<x, b_t(x)\>\le \kk_1-\kk_2|x|^2,\ \ |b_t(x)|\le \kk_2 (1+|x|),\ \  \|\si_{t}(x)\|_{\mathrm{op}}\le \kk_2,\ \ \forall x\in\R^d.\end{equation}
Then for any $p\in (0,\infty)$,  there exists a constant $\kk=\kk(\kk_1,\kk_2,p)\in (0,\infty)$ such that
\beq\label{Z1(2)}\sup_{t\ge 0}\E |X_t^x|^p \le \kk (d^{\ff p 2}+|x|^p),\ \ \forall x\in\R^d.\end{equation}
If moreover
\beq\label{ZZ'(2)} |b_t(x)-b_s(y)|\le K(|x-y|+|x-y|^\aa+|t-s|),\ \ \forall x,y\in \R^d \end{equation}
holds for some constant $K\in (0,\infty),$  then there exists $\kk'=\kk'(\kk_1,\kk_2,K,\eta,\aa,p)\in (0,\infty)$ such that 
\beq\label{Z2(2)}\sup_{t\ge 0}\E |Y_t^x|^p \le \kk' (d^{\ff p 2}+|x|^p),\ \ \forall x\in\R^d.\end{equation}
\end{lem}

\beg{lem}\label{C2.3(2)}
Assume that the conditions in Lemma $\ref{L2.2(2)}$ hold. Then, for any $p>0$ there exists a constant  $\kappa=\kappa(\kk_1,\kk_2,p)\in (0,\infty)$ such that
$$  \E |X_t^x-x|^p \le \kk (d^{\ff p 2}+|x|^{p})(1\land t)^{\ff p 2 \land 1},\ \ x\in\R^d,\ \ t\ge 0.$$
If moreover $\eqref{ZZ'(2)}$ holds,  then there exists $\kk'=\kk'(\kk_1,\kk_2,\eta,K,p,\aa)\in (0,\infty)$ such that
$$\E |Y_t^x-Y_{t_{k-1}}^x|^p \le \kk' d^{\ff p 2}(1+|x|^p)\eta_k^{\ff p 2},\ \ x\in\R^d,\ k\ge 1, \ t\in [t_{k-1},t_k].$$
\end{lem}

\ \newline

\beg{proof}[Proof of Theorem \ref{T1}]   All constants below depend only on $K_1,K_2, K_2',\eta,\aa$ and $p$. \par
(a) In the uniform dissipation case, the well-posedness is well known, {see for instance \cite{Ce01}.
Let $k\ge 0$. By \cite[Theorem 4.1]{Vi09},} we can choose $\F_{t_{k}}$-measurable random variables $\bar X_{t_k}$ and $\bar Y_{t_k}$ such that
 \beq\label{O1} \L(X_{t_k}^x)= \L(\bar X_{t_k}),\ \ \L(Y_{t_k}^x)= \L(\bar Y_{t_k}),\ \ \W_p(\L(X_{t_k}^x),\L(Y_{t_k}^x))^p=\E|\bar X_{t_k}-\bar Y_{t_k}|^p.\end{equation}
 We consider the SDEs for $t\in [t_{k},t_{k+1}]$,
 \beg{align*} &\d \bar X_t= b_t(\bar X_t)\d t+\si_{t}(\bar X_t)\d B_t,\\
 &\d \bar Y_t= b_{t_k}(\bar Y_{t_k})\d t+\si_{t_{k}}(\bar Y_{t_k})\d B_t.\end{align*}
 By the weak uniqueness and the definition of $\W_p$, we obtain
 \beq\label{O2} \W_p\big(\L(X_{t_{k+1}}^x), \L(Y_{t_{k+1}}^x)\big)^p\le \E|\bar X_{t_{k+1}}- \bar Y_{t_{k+1}}|^p.\end{equation}
Let $Z_t= \bar X_t-\bar Y_t.$
 By $(A2)$ and It\^o's formula, we find  a martingale $(M_t)_{t\in [t_k,t_{k+1}]}$ such that 
 \beg{align*}
 & \d |Z_t|^p -\d M_t\\
 &= |Z_t|^{p-2} \Big[p \<Z_t, b_t(\bar X_t)-b_{t_k}(\bar Y_{t_k})\> +\ff p 2 \|\si_{t}(\bar X_t)-\si_{t_{k}}(\bar Y_{t_k})\|_{\mathrm{HS}}^2 \\
 &\qquad \qquad \qquad +\ff 1 2 p(p-2) \ff{|\{\si_{t}(\bar X_t)-\si_{t_{k}}(\bar Y_{t_k})\}^*Z_{t}|^2}{|Z_{t}|^2}  \Big]\d t\\
 &\le |Z_t|^{p-2} \Big[p \<Z_t, b_t(\bar X_t)-b_t(\bar Y_{t})\> +\ff p 2 \|\si_{t}(\bar X_t)-\si_{t}(\bar Y_{t}))\|_{\mathrm{HS}}^2 \\
 &\qquad \qquad \quad  + \ff 1 2 p(p-2)\ff{|\{\si_{t}(\bar X_t)-\si_{t}(\bar Y_{t})\}^*Z_{t}|^2}{|Z_{t}|^2} \\
 &\quad \qquad \qquad + p|Z_t||b_t(\bar Y_t)- b_{t_k}(\bar Y_{t_k})|+ \ff 1 2 p(p-1)\big\{\|\si_{t}(\bar Y_t)-\si_{t_{k}}(\bar Y_{t_k})\|_{\mathrm{HS}}^2\big\}\Big]\d t,\\
 &\le -K_2|Z_t|^p \d t+  pK_1 |Z_t|^{p-1}\Big(|\bar Y_t-\bar Y_{t_{k}}|^\aa+ |\bar Y_t-\bar Y_{t_{k}}|+|t-t_k|\Big)\d t\\
  &\quad +p(p-1) K_1^2 |Z_t|^{p-2}\left[|\bar Y_t- \bar Y_{t_k}|^2+(t-t_{k})^{2} \right]\mathrm{d}t,\ \   t\in [t_k, t_{k+1}].
\end{align*}
By Young's inequality, for any fixed $K_2'\in (0,K_2)$, we  find a constant $c_1>0$ such that
$$\d |Z_t|^p -\d M_t\le - K_2' |Z_t|^p\d t+ c_1\Big(|\bar Y_t-\bar Y_{t_{k}}|^{p\aa}+ |\bar Y_t-\bar Y_{t_{k}}|^p+\eta_{k+1}^{p}\Big)\d t,\ \ t\in [t_k,t_{k+1}],$$
which further implies that, for $t\in [t_k,t_{k+1}]$
\beq\label{O3} \E|Z_t|^p\le \e^{-  K_2' (t-t_k)} \E|Z_{t_k}|^p + c_1 \E\int_{t_k}^t \Big(|\bar Y_s-\bar Y_{t_{k}}|^{\aa p}+ |\bar Y_s-\bar Y_{t_{k}}|^p+\eta_{k+1}^{p}\Big) \d s .\end{equation}
By Lemma \ref{C2.3(2)}, Lemma \ref{L2.2(2)} and $\L(\bar Y_{t_k})=\L(Y_{t_k}^x)$, we can find a constant $c_2>0$ such that
\beq\label{O4} \E|\bar Y_s-\bar Y_{t_{k}}|^p= \E \left[|Y_{s-t_k}^z-z|^p\big{|}z= Y_{t_k}^x\right]  \le c_2 (d\eta_{k+1})^{\ff p 2}(1+|x|^p),\ \ s\in [t_k,t_{k+1}],\end{equation}
so that by Jensen's inequality and \eqref{O3}, we find a constant $c_3>0$ such that
$$  \E|Z_t|^p\le \e^{-  K_2' (t-t_k)} \E|Z_{t_k}|^p + c_3d^{\ff p 2} \eta_{k+1}^{1+\ff{p\aa}2} (1+|x|^p),\ \ t\in [t_k,t_{k+1}].$$
Combining this with \eqref{O1} and  \eqref{O2}, we derive
$$\W_p\big(\L(X_{t_{k+1}}^x),\L(Y_{t_{k+1}}^x)\big)^p\le \e^{- K_2' \eta_{k+1}} \W_p\big(\L_{X_{t_{k}}^x},\L_{Y_{t_{k}}^x}\big)^p+ c_3 d^{\ff p 2}\eta_{k+1}^{1+\ff{p\aa}2} (1+|x|^p),\ \ k\ge 0, x\in\R^d.$$
Iterating in $k$ we prove the desired upper bound estimate in \eqref{MU0}.

(b) Let $\eta_k=\ff\theta k(k\ge 1)$ for some constant $\theta\in (0,\infty)$. Then there exists a constant $c_4\in (0,\infty)$ such that
$$\theta\log\frac{n}{k} -c_4\theta\le t_n-t_k\le \theta \log \ff nk +\theta c_4,\ \ 1\le k\le n.$$
By this and  \eqref{MU0},  when $\theta\ne \ff{\aa p}{2K_2^{\prime}}$, we can find constants $c_5,c_6\in (0,\infty)$ such that
\beg{align*} &\W_p\big(\L(X_{t_n}^x), \L(Y_{t_n}^x)\big)^p
 \le c_4 (1+|x|^p)n^{- \theta   K_2'}\sum_{k=1}^n k^{-1-\ff{\aa p }2   +\theta K_2'}\\
 & \le c_5(1+|x|^p)n^{-\big((\theta K_2') \land \frac{\alpha p}{2}\big)},\ \ n\ge 1, x\in\R^d.\end{align*}
This implies   \eqref{MU} when $\theta\ne\ff{\aa p}{2K_2'}$.  If $\theta=\ff{\aa p}{2K_2'}$,  we may apply the above estimate for $K_2''\in (K_2',K_2)$ replacing $K_2'$,
so that \eqref{MU} holds as well.

(c) Let $(\sigma_{t},b_t) =(\si,b)$ does not depend on $t$. The exponential ergodicity can be proved in a standard way \cite{Ha10}.  By $(A2)$,  synchronuous coupling and It\^o's formula,  for any $x,y\in\R^d$ we can find a martingale $M_t$ such that
$$\d |X_t^x-X_t^y|^p\le -K_2 |X_t^x-X_t^y|^p\d t+\d M_t,\ \ t\ge 0.$$
So, $\E|X_t^x-X_t^y|^p\le |x-y|^p\e^{-K_2t}, t\ge 0.$
By the Markov property, this implies, for any $\mu_1,\mu_2\in \scr P$,

\begin{align}\label{00}
\begin{split}
    \W_p(P_t^*\mu_1,P_t^*\mu_2)^p
&\le  \inf_{\pi\in \C(\mu_{1},\mu_{2})} \int_{\mathbb{R}^{d}\times \mathbb{R}^{d}} \mathbb{E}\left[\left| X^{x}_{t}-X^{y}_{t} \right|^{p} \right] \pi(\mathrm{d}x,\mathrm{d}y)\\
&\le \inf_{\pi\in \C(\mu_{1},\mu_{2})} \int_{\mathbb{R}^{d}\times \mathbb{R}^{d}} \left| x-y \right|^{p}e^{-K_{2}t}  \pi(\mathrm{d}x,\mathrm{d}y)
= \W_p(\mu_1,\mu_2)^p \e^{-K_2 t},
\end{split}
\end{align}
where the first inequality is by $P_t^*\nu=\L(X_t)$ for $X_t$ solving \eqref{SDE} with initial distribution $\nu\in \P$. 
  Moreover, let $L$ be the generator associated with \eqref{SDE} given by
  \begin{align} \label{generator}
      Lf=\langle b, \nabla f \rangle+\frac{1}{2}\langle \sigma \sigma^T, \nabla^{2}f \rangle_{\mathrm{HS}}, \quad \forall f\in C^{2}\left(\mathbb{R}^{d};\mathbb{R} \right).
  \end{align}
  Then $(A2)$ implies that
$$L |\cdot|^{p\lor 2}\le c_1 d^{\ff {p\lor 2} 2}-c_2 |\cdot|^p$$
for some constants $c_1,c_2>0$. By a standard tightness argument, this implies that $P_t^*$ has an invariant probability measure $\mu$ with  $\mu(|\cdot|^p)\le \ff{c_1d^{\ff {p\lor 2} 2}}{c_2}<\infty.$  Combining this with \eqref{00} we conclude that $\mu$ is the unique
  invariant probability measure of $P_t^*$, and
there exists a constant $c_0'\in (0,\infty)$  such that
\beq\label{0}\W_p(\L(X_t^x),\mu)^p=\W_p(P_t^*\dd_x, \mu)^p\le  \mu(|x-\cdot|^p)  \e^{-K_2 t}\le c_0'(1+|x|^p)\e^{-K_2 t} ,\ \ x\in\R^d, t\ge 1.\end{equation}
Combining this together with  (\ref{MU0}) and the triangle inequality, implies the desired upper bound for $\W_p(\L(Y_{t_n}^x),\mu)^p.$
\end{proof}


\begin{thebibliography}{999}

\bibitem{A-J-K} Alfonsi, A., Jourdain, B., and Kohatsu-Higa, A. (2015). Optimal transport bounds between the time-marginals of a multidimensional diffusion and its Euler scheme. {\it Electron. J. Probab. 20} (2015), no. 70, 1–31.

\bibitem{Ce01} Cerrai, S. (Ed.). (2001). Second order PDE’s in finite and infinite dimension: a probabilistic approach. Berlin, Heidelberg: Springer Berlin Heidelberg.


\bibitem{C-S-X} Peng Chen, Qi-Man Shao, and Lihu Xu. A probability approximation framework: Markov process approach. {\it The Annals of Applied Probability, 33}(2), 1619-1659, April 2023.


\bibitem{D17} Dalalyan, A. S. (2017). Theoretical guarantees for approximate sampling from smooth and log-concave densities. {\it Journal of the Royal Statistical Society. Series B (Statistical Methodology)}, 651-676.


\bibitem{D-M17} Durmus, A., and Moulines, E. (2017). Nonasymptotic convergence analysis for the unadjusted Langevin algorithm. {\it The Annals of Applied Probability,27}(3), 1551–1587.

\bibitem{D-M19} Durmus, A., and Moulines, E. (2019). High-dimensional Bayesian inference via the unadjusted Langevin algorithm. {\it Bernoulli 25}(4A), 2019, 2854–2882.

\bibitem{E16} Eberle, A. (2016). Reflection couplings and contraction rates for diffusions. {\it Probability theory and related fields, 166}(3), 851-886.

\bibitem{E-J-H} E, W., Han, J. and Jentzen, A. Deep Learning-Based Numerical Methods for High-Dimensional Parabolic Partial Differential Equations and Backward Stochastic Differential Equations. {\it Commun. Math. Stat}. \textbf{5}, 349–380 (2017). https://doi.org/10.1007/s40304-017-0117-6

\bibitem{Fu98} Fu, W. J. (1998). Penalized regressions: the bridge versus the lasso. {\it Journal of computational and graphical statistics, 7}(3), 397-416.

\bibitem{G-M} Gyöngy, I., and Martínez, T. (2001). On stochastic differential equations with locally unbounded drift. {\it Czechoslovak Mathematical Journal, 51}(4), 763-783.

\bibitem{G-M93} Gelfand, S. B., and Mitter, S. K. (1993). Metropolis-type annealing algorithms for global optimization in $\mathbb{R}^{d}$. {\it SIAM Journal on Control and Optimization, 31}(1), 111-131.


\bibitem{Ha10} Hairer, M. (2010). Convergence of Markov processes. {\it Lecture notes}.

\bibitem{HRW} X. Huang, P. Ren, F.-Y. Wang, \emph{ Probability  distance estimates  between   diffusion processes and applications to singular McKean-Vlasov SDEs,}  arXiv:2304.07562.


\bibitem{L-P02} Lamberton, D., and Pages, G. (2002). Recursive computation of the invariant distribution of a diffusion. {\it Bernoulli}, 367-405.

\bibitem{L-P03}  Lamberton, D., and Pages, G. (2003). Recursive computation of the invariant distribution of a diffusion: the case of a weakly mean reverting drift. {\it Stochastics and dynamics, 3}(04), 435-451.






\bibitem{M-M-N} Mei, S., Montanari, A., and Nguyen, P. M. (2018). A mean field view of the landscape of two-layer neural networks. {\it Proceedings of the National Academy of Sciences, 115}(33), E7665-E7671.

\bibitem{M-F-W-B} Mou, W., Flammarion, N., Wainwright, M. J., and Bartlett, P. L. (2022). Improved bounds for discretization of Langevin diffusions: Near-optimal rates without convexity. {\it Bernoulli, 28}(3), 1577-1601.


\bibitem{N06} Nualart, D. (2006). {\it The Malliavin calculus and related topics} (Vol. 1995, p. 317). Berlin: Springer.

\bibitem{P-P} Pages, G., and Panloup, F. (2023). Unadjusted Langevin algorithm with multiplicative noise: Total variation and Wasserstein bounds. {\it The Annals of Applied Probability, 33}(1), 726-779.

\bibitem{P-W} Priola, E., and Wang, F. Y. (2006). Gradient estimates for diffusion semigroups with singular coefficients. {\it Journal of Functional Analysis, 236}(1), 244-264.

\bibitem{R23} Ren, P. (2023). Singular McKean–Vlasov SDEs: Well-posedness, regularities and Wang’s Harnack inequality. {\it Stochastic Processes and their Applications, 156}, 291-311.

\bibitem{R-T} Roberts, G. O., and Tweedie, R. L. (1996). Exponential convergence of Langevin distributions and their discrete approximations. {\it Bernoulli}, 341-363.

\bibitem{RY13} Revuz, D., and Yor, M. (2013). {\it Continuous martingales and Brownian motion} (Vol. 293). Springer Science and Business Media.

\bibitem{Vi09} Villani, C. (2009). {\it Optimal transport: old and new} (Vol. 338, p. 23). Berlin: springer.




\bibitem{W23} F.-Y. Wang, \emph{Exponential Ergodicity for singular reflecting  McKean-Vlasov SDEs,}  Stoch. Proc. Apll. 160(2023), 265--293.

\bibitem{X-Z} Xie, L., and Zhang, X. (2016). Sobolev differentiable flows of SDEs with local Sobolev and super-linear growth coefficients. {\it The Annals of Probability, 44}(6), 3661-3687.

\bibitem{XZ} L. Xie, X. Zhang, \emph{ Ergodicity of stochastic differential equations with jumps and singular coefficients,}  Ann. Inst. Henri Poincar\'e Probab. Stat. 56(2020),  175--229.

\bibitem{Z05} Zhang, X. (2005). Strong solutions of SDES with singular drift and Sobolev diffusion coefficients. {\it Stochastic Processes and their Applications, 115}(11), 1805-1818.

\bibitem{Zhang11} X. Zhang, \emph{Stochastic homeomorphism flows of SDEs with singular drifts and Sobolev diffusion coefficients, }  Electr. J. Probab.
 16(2011), 1096--1116.

























\end{thebibliography}
\end{document}